\documentclass[12pt]{article}
	\newcommand{\mainTitle}{A parameterized generalization of the sum formula for quadruple zeta values}

	\newcommand{\authorName}{MACHIDE, Tomoya}
	\newcommand{\organizationName}{Kinki University}
	\newcommand{\departmentName}{Interdisciplinary Graduate School of Science and Engineering}
	\newcommand{\majorName}{Research Center for Quantum Computing}
	\newcommand{\placeAddress}{3-4-1 Kowakae, Higashi-Osaka, Osaka 577-8502, Japan}
	\newcommand{\emailAddress}{E-mail: machide.t@gmail.com}
	
	\newcommand{\sectOne}{Introduction}
	\newcommand{\sectTwo}{Parameterized sums of multiple polylogarithms}
	\newcommand{\sSectTwoO}{An identity derived from harmonic relations}
	\newcommand{\sSectTwoT}{Identities derived from shuffle relations}
	\newcommand{\sSectTwoTh}{Asymptotic properties}
	\newcommand{\sectThree}{Proof of \refThm{1_Thm1}}
	\newcommand{\sectFour}{Proof of \refThm{1_Thm2}}

\usepackage{geometry}               
\usepackage{ascmac}
\usepackage{amsmath}
\usepackage{amssymb}
\usepackage{amsthm}
\usepackage[hypertex, colorlinks=true, bookmarks=true]{hyperref}

	\DeclareMathOperator*{\OPlus}{\bigoplus}
	\newcommand{\nbk}[3]{#1#3#2}		
	\newcommand{\bgbk}[3]{\bigl{#1}#3\bigr{#2}}	
	\newcommand{\Bgbk}[3]{\Bigl{#1}#3\Bigr{#2}}			
	\newcommand{\bggbk}[3]{\biggl{#1}#3\biggr{#2}}			
	\newcommand{\Bggbk}[3]{\Biggl{#1}#3\Biggr{#2}}
	\newcommand{\autobk}[3]{\left#1#3\right#2}
	\newcommand{\mcbk}[4][?]{\ifx n#1\nbk{#2}{#3}{#4}\else\ifx b#1\bgbk{#2}{#3}{#4}\else\ifx B#1\Bgbk{#2}{#3}{#4}\else\ifx g#1\bggbk{#2}{#3}{#4}\else\ifx G#1\Bggbk{#2}{#3}{#4}\else\ifx a#1\autobk{#2}{#3}{#4}\else#4\fi\fi\fi\fi\fi\fi}
	\newcommand{\nsgsb}[1]{#1}		
	\newcommand{\bgsgsb}[1]{\big{#1}}	
	\newcommand{\Bgsgsb}[1]{\Big{#1}}			
	\newcommand{\bggsgsb}[1]{\bigg{#1}}			
	\newcommand{\Bggsgsb}[1]{\Bigg{#1}}
	\newcommand{\mcsgsb}[2][?]{\ifx n#1\nsgsb{#2}\else\ifx b#1\bgsgsb{#2}\else\ifx B#1\Bgsgsb{#2}\else\ifx g#1\bggsgsb{#2}\else\ifx G#1\Bggsgsb{#2}\else#2\fi\fi\fi\fi\fi}
	\newcommand{\myEqSpace}{\,}
	\newlength{\myEqSpaceLen}
	\newcommand{\bkR}[2][a]{\mcbk[#1]{(}{)}{#2}}
	\newcommand{\bkS}[2][a]{\mcbk[#1]{[}{]}{#2}}
	\newcommand{\bkB}[2][a]{\mcbk[#1]{\{}{\}}{#2}}
	\newcommand{\bkA}[2][a]{\mcbk[#1]{\langle}{\rangle}{#2}}

	\newcommand{\tmR}[4][?]{\bkR[#1]{#2} \ifx 0#4 \else \ifx 1#4 #3 \else {#3}^{#4} \fi \fi}	
	\newcommand{\tmS}[4][?]{\bkS[#1]{#2} \ifx 0#4 \else \ifx 1#4 #3 \else {#3}^{#4} \fi \fi}
	\newcommand{\tmM}[4][?]{\bkB[#1]{#2} \ifx 0#4 \else \ifx 1#4 #3 \else {#3}^{#4} \fi \fi}
		
	\newcommand{\setN}{\mathbb{N}}

	\newcommand{\setQ}{\mathbb{Q}}

	\newcommand{\setC}{\mathbb{C}}
	
	\newcommand{\gpSym}[2][?]{\ifx g#1\mathfrak{S}_{#2} \else S_{#2} \fi}
	\newcommand{\gpAlt}[2][?]{\ifx g#1 \mathfrak{A}_{#2} \else A_{#2} \fi }
	\newcommand{\gpKleinF}[1][?]{\ifx g#1 \mathfrak{V} \else V \fi }
	\newcommand{\gpCyc}[2][n]{\bkA[#1]{#2}}

	\newcommand{\gpu}{e}
	
	\newcommand{\vPack}[1][10]{\vspace{-#1pt}}
	
	\newcommand{\lnA}[1][]{&  &}
	\newcommand{\lnP}[1]{\myEqSpace#1\myEqSpace}
	\newcommand{\lnAP}[2][]{& #2 &}
	\newcommand{\lnAH}[1][\nonumber]{#1\\ & &}
	\newcommand{\lnAHP}[2][\nonumber]{#1 \\ & #2 &}
	\newcommand{\lnAHs}[2][\nonumber]{#1 \\ & &\hspace{#2pt}}

		\newcommand{\slnAHs}[1]{\\ \hspace{#1pt}}
	\newcommand{\refH}[2]{#1\ref{#2}}
	
	\newcommand{\refEq}[1]{(\ref{#1})}	
		
	\newcommand{\refSect}[2][?]{\refH{\S}{#2}}
	
	\newcommand{\alp}{\alpha}

	\newcommand{\sig}{\sigma}

	\newcommand{\mVert}[1][n]{{\,\mcsgsb[#1]{\vert}\,}}
	
	\newcommand{\opF}[3][?]{\ifx s#1#2/#3\else\ifx b#1(#2)/(#3)\else\ifx d#1\dfrac{#2}{#3}\else\frac{#2}{#3}\fi\fi\fi}

	\newcommand{\pw}[3][?]{\ifx!#3{#2}^{#3}\else#2^{#3}\fi}
	\newcommand{\id}[3][?]{#2_{#3}}
	\newcommand{\ip}[4][?]{{#2}_{#3}^{#4}}
	\newcommand{\pwR}[3][a]{\ifx!#1{\bkR[#1]{#2}}^{#3}\else\bkR[#1]{#2}^{#3}\fi}
	\newcommand{\pwB}[3][a]{\ifx!#1{\bkB[#1]{#2}}^{#3}\else\bkB[#1]{#2}^{#3}\fi}
	\newcommand{\pwS}[3][a]{\ifx!#1{\bkS[#1]{#2}}^{#3}\else\bkS[#1]{#2}^{#3}\fi}
					
	\newcommand{\nFc}[3][n]{#2\bkR[#1]{#3}}
	
	\newcommand{\idFc}[4][n]{\id{#2}{#3}\bkR[#1]{#4}}
	\newcommand{\pwFc}[4][n]{\pw{#2}{#3}\bkR[#1]{#4}}
	\newcommand{\ipFc}[5][n]{\ip{#2}{#3}{#4}\bkR[#1]{#5}}
		\newcommand{\Fc}{\nFc}
	
	\newcommand{\tpT}[3][a]{ {#2}\atop \bkR[#1]{#3} }

	\newcommand{\nSm}[2][?]{\ifx l#1 \sum\limits_{#2} \else \sum_{#2} \fi}
	\newcommand{\nSmT}[3][?]{\ifx l#1 \sum\limits_{#2}^{#3} \else \sum_{#2}^{#3} \fi}	
	\newcommand{\pSm}[2][?]{\ifx t#1 \sum_{#2}^{\prime} \else \sideset{}{^\prime}\sum_{#2} \fi}
	\newcommand{\pSmT}[3][?]{\ifx t#1 \sum_{#2}^{\prime#3} \else \sideset{}{^\prime}\sum_{#2}^{#3} \fi}	
	\newcommand{\pSmN}[1][?]{\ifx t#1 \sum^{\prime} \else \sideset{}{^\prime}\sum \fi}
	\newcommand{\dSm}[2][?]{\ifx t#1 \sum_{#2}^{\dagger} \else \sideset{}{^\dagger}\sum_{#2} \fi}
	\newcommand{\dSmT}[3][?]{\ifx t#1 \sum_{#2}^{\dagger#3} \else \sideset{}{^\dagger}\sum_{#2}^{#3} \fi}	
	\newcommand{\dSmN}[1][?]{\ifx t#1 \sum^{\dagger} \else \sideset{}{^\dagger}\sum \fi}
	\newcommand{\tpTSm}[3][?]{\nSm[#1]{\tpT{#2}{#3}}}

		\newcommand{\Sm}{\nSm}
		\newcommand{\SmT}{\nSmT}
		\newcommand{\tpSm}{\tpTSm}
		
	\newcommand{\nPd}[2][?]{\ifx l#1 \prod\limits_{#2} \else \prod_{#2} \fi}
	\newcommand{\nPdT}[3][?]{\ifx l#1 \prod\limits_{#2}^{#3} \else \prod_{#2}^{#3} \fi}

	\newcommand{\nOPs}[2][?]{\ifx l#1 \OPlus\limits_{#2} \else \OPlus_{#2} \fi}
	\newcommand{\nOPsT}[3][?]{\ifx l#1 \OPlus\limits_{#2}^{#3} \else \OPlus_{#2}^{#3} \fi}	
	\newcommand{\pOPs}[2][?]{\ifx t#1 \OPlus_{#2}^{\prime} \else \sideset{}{^\prime}\OPlus_{#2} \fi}
	\newcommand{\pOPsT}[3][?]{\ifx t#1 \OPlus_{#2}^{\prime#3} \else \sideset{}{^\prime}\OPlus_{#2}^{#3} \fi}

	\newcommand{\nIs}[2][?]{\ifx l#1 \bigcap\limits_{#2} \else \bigcap_{#2} \fi}
	\newcommand{\nIsT}[3][?]{\ifx l#1 \bigcap\limits_{#2}^{#3} \else \bigcap_{#2}^{#3} \fi}	
	\newcommand{\pIs}[2][?]{\ifx t#1 \bigcap_{#2}^{\prime} \else \sideset{}{^\prime}\bigcap_{#2} \fi}
	\newcommand{\pIsT}[3][?]{\ifx t#1 \bigcap_{#2}^{\prime#3} \else \sideset{}{^\prime}\bigcap_{#2}^{#3} \fi}

	\newcommand{\nUn}[2][?]{\ifx l#1 \bigcup\limits_{#2} \else \bigcup_{#2} \fi}
	\newcommand{\nUnT}[3][?]{\ifx l#1 \bigcup\limits_{#2}^{#3} \else \bigcup_{#2}^{#3} \fi}	
	\newcommand{\pUn}[2][?]{\ifx t#1 \bigcup_{#2}^{\prime} \else \sideset{}{^\prime}\bigcup_{#2} \fi}
	\newcommand{\pUnT}[3][?]{\ifx t#1 \bigcup_{#2}^{\prime#3} \else \sideset{}{^\prime}\bigcup_{#2}^{#3} \fi}

	\newcommand{\nLm}[2][?]{\ifx l#1 \lim\limits_{#2} \else \lim_{#2} \fi}
		\newcommand{\Lm}{\nLm}
	\newcommand{\nMax}[2][a]{ \max\bkB[#1]{#2} }
	
		\newcommand{\Max}{\nMax}
		
	\newcommand{\glcondEnvLineHead}[1]{ \ifx*#1 \begin{eqnarray*} \else \begin{eqnarray}  \label{#1} \fi }
	\newcommand{\glcondEnvLineTail}[1]{ \ifx*#1 \end{eqnarray*} \else \end{eqnarray} \fi }
	\newcommand{\glcondDis}[1]{\ifx d#1 \displaystyle \fi}

		\newcommand{\envHLineT}[3][*]{ \glcondEnvLineHead{#1} #2&=&#3\glcondEnvLineTail{#1} }
		\newcommand{\envHLineTDef}[3][*]{ \glcondEnvLineHead{#1} #2&:=&#3\glcondEnvLineTail{#1} }
			\newcommand{\envHLine}{\envHLineT}
			\newcommand{\envHLineDef}{\envHLineTDef}
		
		\newcommand{\envHLineThDef}[4][*]{ \glcondEnvLineHead{#1} #2&:=&#3\\&=&#4\nonumber \glcondEnvLineTail{#1} }
		\newcommand{\envHLineF}[5][*]{ \glcondEnvLineHead{#1} #2&=&#3\\&=&#4\nonumber \\&=&#5\nonumber \glcondEnvLineTail{#1} }
		
		\newcommand{\envHLineCSNme}[7][*]{\begin{eqnarray} #2&=&#3\\#4&=&#5\\#6&=&#7\end{eqnarray}}
		
		\newcommand{\envPLine}[2][*]{\glcondEnvLineHead{#1} #2\glcondEnvLineTail{#1}}
		
		\newcommand{\envSLineO}[2][*]{\renewcommand{\arraystretch}{1.4} \begin{array}{l}#2\end{array} \renewcommand{\arraystretch}{1}}
			\newcommand{\envSLine}{\envSLineO}

		\newcommand{\envMO}[2][*]{$\ifx d#1 \displaystyle \fi#2$}
		\newcommand{\envMT}[3][*]{$\ifx d#1 \displaystyle \fi#2=#3$}
		\newcommand{\envMTDef}[3][*]{$\ifx d#1 \displaystyle \fi#2:=#3$}
			\newcommand{\envM}{\envMT}
			\newcommand{\envMDef}{\envMTDef}
		\newcommand{\envMTh}[4][*]{$\ifx d#1 \displaystyle \fi#2=#3=#4$}
		\newcommand{\envMF}[5][*]{$\ifx d#1 \displaystyle \fi#2=#3=#4=#5$}
		
	\newcommand{\envMLineT}[3][*]{ \ifx*#1 \begin{multline*} #2\lnP{=}#3\end{multline*} \else \begin{multline} \label{#1} #2\lnP{=}#3\end{multline} \fi }
	\newcommand{\envMLineTDef}[3][*]{ \ifx*#1 \begin{multline*} #2\lnP{:=}#3\end{multline*} \else \begin{multline} \label{#1} #2\lnP{:=}#3\end{multline} \fi }
		\newcommand{\envMLine}{\envMLineT}

	\newcommand{\lcparaCase}{\vspace{3pt}}

	\newcommand{\abs}[1]{\left | #1 \right |  }		

	\newcommand{\sbLandau}[2][n]{\Fc[#1]{O}{#2}}

	\newcommand{\lgg}[2][?]{\Fc[#1]{\log}{#2}}
	
	\newcommand{\fcZeta}[2][n]{\Fc[#1]{\zeta}{#2}}
		\newcommand{\fcZ}{\fcZeta}		
	\newcommand{\lgP}[3][n]{\idFc[#1]{Li}{#2}{#3}}
	
		\newcommand{\envMTPt}[4][*]{$\ifx d#1 \displaystyle \fi#3#2#4$}
			
		\newcommand{\envMThPt}[5][*]{$\ifx d#1 \displaystyle \fi#3#2#4#2#5$}
		\newcommand{\envMFPt}[6][*]{$\ifx d#1 \displaystyle \fi#3#2#4#2#5#2#6$}

	\newcommand{\envMLineTPt}[4][*]{ \ifx*#1 \begin{multline*} #3\lnP{#2}#4\end{multline*} \else \begin{multline} \label{#1} #3\lnP{#2}#4\end{multline} \fi }

	\newcommand{\cTxT}[2]{\textcolor{#1}{#2}}
	
		\newcommand{\cTx}{\cTxT}
	\newcommand{\alTx}[1]{\cTx{red}{#1}}
	\newcommand{\rmTx}[1]{\cTx{blue}{#1}}
	\newcommand{\ntTx}[1]{\cTx{Green}{#1}}
	\newcommand{\cmoTx}[1]{\cTx{Gray}{#1}}
	\newcommand{\sTx}[2][n]{ \ifx t#1{\tiny #2} \else \ifx s#1{\scriptsize #2} \else \ifx f#1{\footnotesize #2} \else \ifx S#1{\small #2} \else \ifx n#1{#2} \else \ifx l#1{\large #2} \else \ifx L#1{\Large #2} \else \ifx R#1{\LARGE #2} \else \ifx h#1{\huge #2} \else \ifx H#1{\Huge #2} \else \ifx ?#1 #2 \else #2 \fi\fi\fi\fi\fi\fi\fi\fi\fi\fi\fi }
	\newcommand{\bfTx}[1]{{\bf#1}}

		\newcommand{\envHLineTCl}[3][a]{ \ifx a#1\alTx{\envHLineT{#2}{#3}} \else \ifx r#1 \rmTx{\envHLineT{#2}{#3}} \else\ifx n#1 \ntTx{\envHLineT{#2}{#3}}\else\ifx c#1 \cmoTx{\envHLineT{#2}{#3}} \else \text{[argument error]} \fi\fi\fi\fi \vPack[18] }

		\newcommand{\envHLineCSClPart}[8][?]{\ifx*#1 \begin{eqnarray*} \else \begin{eqnarray}  \label{#1}  \fi \alTx{#3}&#2&\alTx{#4}\\#5&#2&#6\nonumber\\\alTx{#7}&#2&\alTx{#8}\nonumber\glcondEnvLineTail{*}}
		
			\newcommand{\HLineCTCl}[3][?]{\alTx{#2}&=&\alTx{#3}\nonumber \ifx#1p \\ \fi}
			\newcommand{\HLineCTClDef}[3][?]{\alTx{#2}&:=&\alTx{#3}\nonumber \ifx#1p \\ \fi}
			\newcommand{\HLineCFCl}[5][?]{\alTx{#2}&=&\alTx{#3}\nonumber\\#4&=&#5\nonumber \ifx#1p \\ \fi}
			\newcommand{\HLineCFClDef}[5][?]{\alTx{#2}&:=&\alTx{#3}\nonumber\\#4&:=&#5\nonumber \ifx#1p \\ \fi}
			\newcommand{\HLineCSCl}[7][?]{\alTx{#2}&=&\alTx{#3}\nonumber\\#4&=&#5\nonumber\\\alTx{#6}&=&\alTx{#7}\nonumber\ifx#1p \\ \fi}
			\newcommand{\HLineCSClDef}[7][?]{\alTx{#2}&:=&\alTx{#3}\nonumber\\#4&:=&#5\nonumber\\\alTx{#6}&:=&\alTx{#7}\nonumber\ifx#1p \\ \fi}									
			\newcommand{\HLineCECl}[9][?]{\alTx{#2}&=&\alTx{#3}\nonumber\\#4&=&#5\nonumber\\\alTx{#6}&=&\alTx{#7}\nonumber\\#8&=&#9\nonumber\ifx#1p \\ \fi}
			\newcommand{\HLineCEClDef}[9][?]{\alTx{#2}&:=&\alTx{#3}\nonumber\\#4&:=&#5\nonumber\\\alTx{#6}&:=&\alTx{#7}\nonumber\\#8&:=&#9\nonumber\ifx#1p \\ \fi}
		
\theoremstyle{plain}
\newtheorem{theorem}{THEOREM}[section]
\newtheorem{proposition}[theorem]{PROPOSITION}
\newtheorem{lemma}[theorem]{LEMMA}

\theoremstyle{definition}

\theoremstyle{remark}
\newtheorem{remark}[theorem]{REMARK}
\allowdisplaybreaks[4]
\numberwithin{equation}{section}

	\newcommand{\lccondBibitem}[3][]{ \if ?#2 \bibitem{#3} \else \bibitem[#2]{#3} \fi}
	\newcommand{\refPaper}[8][?]{
			\lccondBibitem{#1}{#2}
				#3,			
				\emph{#4}, 	
				#5\ 			
				{\bf #6},		
				#7,			
				#8.			
		}
	\newcommand{\refPreprint}[6][?]{
			\lccondBibitem{#1}{#2}
				#3,			
				\emph{#4}, 	
				preprint; #5,			
				#6.			
		}
	\newcommand{\refPaperRep}[9][?]{
			\lccondBibitem{#1}{#2}
				#3,			
				\emph{#4}, 	
				#5\ 			
				{\bf #6},		
				#7,			
				#8			
				; reprinted in #9	
		}

	\newcommand{\refPaperAlm}[5][?]{
			\lccondBibitem{#1}{#2}
				#3,	 		
				\emph{#4}, 	
				#5		
		}

	\newcommand{\pcstSpaceForRef}{\ }
	\newcommand{\refThm}[2][?]{\ifx?#1\refH{Theorem\pcstSpaceForRef}{#2}\else\ifx s#1\refH{Theorems\pcstSpaceForRef}{#2}\else{[argument error]}\fi\fi}
	\newcommand{\refProp}[2][?]{\ifx?#1\refH{Proposition\pcstSpaceForRef}{#2}\else\ifx s#1\refH{Propositions\pcstSpaceForRef}{#2}\else{[argument error]}\fi\fi}
	\newcommand{\refLem}[2][?]{\ifx?#1\refH{Lemma\pcstSpaceForRef}{#2}\else\ifx s#1\refH{Lemmas\pcstSpaceForRef}{#2}\else{[argument error]}\fi\fi}
	\newcommand{\refCor}[2][?]{\ifx?#1\refH{Corollary\pcstSpaceForRef}{#2}\else\ifx s#1\refH{Corollaries\pcstSpaceForRef}{#2}\else{[argument error]}\fi\fi}
	\newcommand{\refDef}[2][?]{\ifx?#1\refH{Definition\pcstSpaceForRef}{#2}\else\ifx s#1\refH{Definitions\pcstSpaceForRef}{#2}\else{[argument error]}\fi\fi}
	\newcommand{\refRem}[2][?]{\ifx?#1\refH{Remark\pcstSpaceForRef}{#2}\else\ifx s#1\refH{Remarks\pcstSpaceForRef}{#2}\else{[argument error]}\fi\fi}
	\newcommand{\refTab}[2][?]{\ifx?#1\refH{Table\pcstSpaceForRef}{#2}\else\ifx s#1\refH{Tables\pcstSpaceForRef}{#2}\else{[argument error]}\fi\fi}
	\newcommand{\rTx}[2][0.5]{ \raise#1ex\hbox{$ \displaystyle#2$} }

	\newcommand{\glcondEnvLineTailPd}[1]{.\ifx*#1 \end{eqnarray*} \else \end{eqnarray} \fi }
	\newcommand{\glcondEnvLineTailCm}[1]{,\ifx*#1 \end{eqnarray*} \else \end{eqnarray} \fi }

	\newcommand{\envProof}[2][?]{ \par\mbox{}\vspace{-5pt}\\ \ifx?#1\emph{Proof.}\else\emph{Proof of #1.}\fi \ #2 \hfill $\Box$\\ \par}
	
		\newcommand{\envLineTPd}[3][*]{ \glcondEnvLineHead{#1} & &#2\\&=&#3\nonumber \glcondEnvLineTailPd{#1} }
		
		\newcommand{\envLineTCm}[3][*]{ \glcondEnvLineHead{#1} & &#2\\&=&#3\nonumber \glcondEnvLineTailCm{#1} }
		
			\newcommand{\envLinePd}{\envLineTPd}
			
			\newcommand{\envLineCm}{\envLineTCm}
			
		\newcommand{\envLineThPd}[4][*]{ \glcondEnvLineHead{#1} & &#2\\&=&#3\nonumber \\&=&#4\nonumber \glcondEnvLineTailPd{#1} }
		\newcommand{\envLineThCm}[4][*]{ \glcondEnvLineHead{#1} & &#2\\&=&#3\nonumber \\&=&#4\nonumber \glcondEnvLineTailCm{#1} }

		\newcommand{\envHLineTPd}[3][*]{ \glcondEnvLineHead{#1} #2&=&#3\glcondEnvLineTailPd{#1} }
		\newcommand{\envHLineTDefPd}[3][*]{ \glcondEnvLineHead{#1} #2&:=&#3\glcondEnvLineTailPd{#1} }
		\newcommand{\envHLineTCm}[3][*]{ \glcondEnvLineHead{#1} #2&=&#3\glcondEnvLineTailCm{#1} }
		
			\newcommand{\envHLinePd}{\envHLineTPd}
			\newcommand{\envHLineDefPd}{\envHLineTDefPd}
			\newcommand{\envHLineCm}{\envHLineTCm}
			
		\newcommand{\envHLineFCm}[5][*]{ \glcondEnvLineHead{#1} #2&=&#3\\&=&#4\nonumber \\&=&#5\nonumber \glcondEnvLineTailCm{#1} }

		\newcommand{\envHLineCFNmePd}[5][?]{\begin{eqnarray} #2&=&#3,\\#4&=&#5 \glcondEnvLineTailPd{?} }
		\newcommand{\envHLineCFDefNmePd}[5][?]{\begin{eqnarray} #2&:=&#3,\\#4&:=&#5 \glcondEnvLineTailPd{?} }
		\newcommand{\envHLineCFCmNme}[5][?]{\begin{eqnarray} #2&=&#3,\\#4&=&#5 \glcondEnvLineTailCm{?} }
		\newcommand{\envHLineCFCmDefNme}[5][?]{\begin{eqnarray} #2&:=&#3,\\#4&:=&#5 \glcondEnvLineTailCm{?} }
		\newcommand{\envHLineCSPd}[7][*]{\glcondEnvLineHead{#1} #2&=&#3,\\#4&=&#5,\nonumber\\#6&=&#7\nonumber\glcondEnvLineTailPd{#1}}
		
		\newcommand{\envHLineCSCm}[7][*]{\glcondEnvLineHead{#1} #2&=&#3,\\#4&=&#5,\nonumber\\#6&=&#7\nonumber\glcondEnvLineTailCm{#1}}
		\newcommand{\envHLineCSCmDef}[7][*]{\glcondEnvLineHead{#1} #2&:=&#3,\\#4&:=&#5,\nonumber\\#6&:=&#7\nonumber\glcondEnvLineTailCm{#1}}	
		\newcommand{\envHLineCSNmePd}[7][*]{\begin{eqnarray} #2&=&#3,\\#4&=&#5,\\#6&=&#7\glcondEnvLineTailPd{?}}
		\newcommand{\envHLineCSDefNmePd}[7][*]{\begin{eqnarray} #2&:=&#3,\\#4&:=&#5,\\#6&:=&#7\glcondEnvLineTailPd{?}}
		\newcommand{\envHLineCSCmNme}[7][*]{\begin{eqnarray} #2&=&#3,\\#4&=&#5,\\#6&=&#7\glcondEnvLineTailCm{?}}
		\newcommand{\envHLineCSCmDefNme}[7][*]{\begin{eqnarray} #2&:=&#3,\\#4&:=&#5,\\#6&:=&#7\glcondEnvLineTailCm{?}}
		\newcommand{\envHLineCEPd}[9][*]{\glcondEnvLineHead{#1} #2&=&#3,\\#4&=&#5,\nonumber\\#6&=&#7,\nonumber\\#8&=&#9\nonumber\glcondEnvLineTailPd{#1}}

		\newcommand{\envHLineCENmePd}[9][*]{\begin{eqnarray} #2&=&#3,\\#4&=&#5,\\#6&=&#7,\\#8&=&#9\glcondEnvLineTailPd{?}}
		\newcommand{\envHLineCEDefNmePd}[9][*]{\begin{eqnarray} #2&:=&#3,\\#4&:=&#5,\\#6&:=&#7,\\#8&:=&#9\glcondEnvLineTailPd{?}}
		\newcommand{\envHLineCECmNme}[9][*]{\begin{eqnarray} #2&=&#3,\\#4&=&#5,\\#6&=&#7,\\#8&=&#9\glcondEnvLineTailCm{?}}
		\newcommand{\envHLineCECmDefNme}[9][*]{\begin{eqnarray} #2&:=&#3,\\#4&:=&#5,\\#6&:=&#7,\\#8&:=&#9\glcondEnvLineTailCm{?}}
		\newcommand{\envPLinePd}[2][*]{\glcondEnvLineHead{#1} #2\glcondEnvLineTailPd{#1}}
		\newcommand{\envPLineCm}[2][*]{\glcondEnvLineHead{#1} #2\glcondEnvLineTailCm{#1}}

		\newcommand{\envOTLineCm}[4][*]{\glcondEnvLineHead{#1} #2\lnAP{=}#3\lnP{=}#4,\glcondEnvLineTail{#1}}

		\newcommand{\envMOCm}[2][*]{$\ifx d#1 \displaystyle \fi#2$,}
		\newcommand{\envMOPd}[2][*]{$\ifx d#1 \displaystyle \fi#2$.}
		
		\newcommand{\envMTCm}[3][*]{$\ifx d#1 \displaystyle \fi#2=#3$,}
		\newcommand{\envMTPd}[3][*]{$\ifx d#1 \displaystyle \fi#2=#3$.}
		\newcommand{\envMTCmDef}[3][*]{$\ifx d#1 \displaystyle \fi#2:=#3$,}
		\newcommand{\envMTDefPd}[3][*]{$\ifx d#1 \displaystyle \fi#2:=#3$.}
			\newcommand{\envMCm}{\envMTCm}
			\newcommand{\envMPd}{\envMTPd}
			\newcommand{\envMCmDef}{\envMTCmDef}
			
		\newcommand{\envMThCm}[4][*]{$\ifx d#1 \displaystyle \fi#2=#3=#4$,}
		\newcommand{\envMThPd}[4][*]{$\ifx d#1 \displaystyle \fi#2=#3=#4$.}

		\newcommand{\envHLineCFCmNm}[5][*]{ \begin{equation}\begin{split} \ifx*#1 \text{[ERROR;need label name]} \else \label{#1} \fi #2&\lnP{=}#3,\\#4&\lnP{=}#5, \end{split}\end{equation} }
		\newcommand{\envHLineCFNm}[5][*]{ \begin{equation}\begin{split} \ifx*#1 \text{[ERROR;need label name]} \else \label{#1} \fi #2&\lnP{=}#3\\#4&\lnP{=}#5, \end{split}\end{equation} }
		\newcommand{\envHLineCFNmPd}[5][*]{ \begin{equation}\begin{split} \ifx*#1 \text{[ERROR;need label name]} \else \label{#1} \fi #2&\lnP{=}#3,\\#4&\lnP{=}#5. \end{split}\end{equation} }
		\newcommand{\envHLineCFCmDefNm}[5][*]{ \begin{equation}\begin{split} \ifx*#1 \text{[ERROR;need label name]} \else \label{#1} \fi #2&\lnP{:=}#3,\\#4&\lnP{:=}#5, \end{split}\end{equation} }
		\newcommand{\envHLineCFDefNm}[5][*]{ \begin{equation}\begin{split} \ifx*#1 \text{[ERROR;need label name]} \else \label{#1} \fi #2&\lnP{:=}#3\\#4&\lnP{:=}#5, \end{split}\end{equation} }
		\newcommand{\envHLineCFDefNmPd}[5][*]{ \begin{equation}\begin{split} \ifx*#1 \text{[ERROR;need label name]} \else \label{#1} \fi #2&\lnP{:=}#3,\\#4&\lnP{:=}#5. \end{split}\end{equation} }
		\newcommand{\envHLineCSCmNm}[7][*]{ \begin{equation}\begin{split} \ifx*#1 \text{[ERROR;need label name]} \else \label{#1} \fi #2&\lnP{=}#3,\\#4&\lnP{=}#5,\\#6&\lnP{=}#7 \end{split}\end{equation} }
		\newcommand{\envHLineCSNm}[7][*]{ \begin{equation}\begin{split} \ifx*#1 \text{[ERROR;need label name]} \else \label{#1} \fi #2&\lnP{=}#3\\#4&\lnP{=}#5\\#6&\lnP{=}#7 \end{split}\end{equation} }
		\newcommand{\envHLineCSNmPd}[7][*]{ \begin{equation}\begin{split} \ifx*#1 \text{[ERROR;need label name]} \else \label{#1} \fi #2&\lnP{=}#3,\\#4&\lnP{=}#5,\\#6&\lnP{=}#7. \end{split}\end{equation} }
		\newcommand{\envHLineCSCmDefNm}[7][*]{ \begin{equation}\begin{split} \ifx*#1 \text{[ERROR;need label name]} \else \label{#1} \fi #2&\lnP{:=}#3,\\#4&\lnP{:=}#5,\\#6&\lnP{:=}#7 \end{split}\end{equation} }
		\newcommand{\envHLineCSDefNm}[7][*]{ \begin{equation}\begin{split} \ifx*#1 \text{[ERROR;need label name]} \else \label{#1} \fi #2&\lnP{:=}#3\\#4&\lnP{:=}#5\\#6&\lnP{:=}#7 \end{split}\end{equation} }
		\newcommand{\envHLineCSDefNmPd}[7][*]{ \begin{equation}\begin{split} \ifx*#1 \text{[ERROR;need label name]} \else \label{#1} \fi #2&\lnP{:=}#3,\\#4&\lnP{:=}#5,\\#6&\lnP{:=}#7. \end{split}\end{equation} }
		\newcommand{\envHLineCECmNm}[9][*]{ \begin{equation}\begin{split} \ifx*#1 \text{[ERROR;need label name]} \else \label{#1} \fi #2&\lnP{=}#3,\\#4&\lnP{=}#5,\\#6&\lnP{=}#7,\\#8&\lnP{=}#9,  \end{split}\end{equation} }
		\newcommand{\envHLineCENm}[9][*]{ \begin{equation}\begin{split} \ifx*#1 \text{[ERROR;need label name]} \else \label{#1} \fi #2&\lnP{=}#3\\#4&\lnP{=}#5\\#6&\lnP{=}#7\\#8&\lnP{=}#9  \end{split}\end{equation} }
		\newcommand{\envHLineCENmPd}[9][*]{ \begin{equation}\begin{split} \ifx*#1 \text{[ERROR;need label name]} \else \label{#1} \fi #2&\lnP{=}#3,\\#4&\lnP{=}#5,\\#6&\lnP{=}#7,\\#8&\lnP{=}#9.  \end{split}\end{equation} }
		\newcommand{\envHLineCECmDefNm}[9][*]{ \begin{equation}\begin{split} \ifx*#1 \text{[ERROR;need label name]} \else \label{#1} \fi #2&\lnP{:=}#3,\\#4&\lnP{:=}#5,\\#6&\lnP{:=}#7,\\#8&\lnP{:=}#9,  \end{split}\end{equation} }
		\newcommand{\envHLineCEDefNm}[9][*]{ \begin{equation}\begin{split} \ifx*#1 \text{[ERROR;need label name]} \else \label{#1} \fi #2&\lnP{:=}#3\\#4&\lnP{:=}#5\\#6&\lnP{:=}#7\\#8&\lnP{:=}#9  \end{split}\end{equation} }
		\newcommand{\envHLineCEDefNmPd}[9][*]{ \begin{equation}\begin{split} \ifx*#1 \text{[ERROR;need label name]} \else \label{#1} \fi #2&\lnP{:=}#3,\\#4&\lnP{:=}#5,\\#6&\lnP{:=}#7,\\#8&\lnP{:=}#9.  \end{split}\end{equation} }
	\newcommand{\envMLineTPd}[3][*]{ \ifx*#1 \begin{multline*} #2\lnP{=}#3.\end{multline*} \else \begin{multline} \label{#1} #2\lnP{=}#3.\end{multline} \fi }
	\newcommand{\envMLineTCm}[3][*]{ \ifx*#1 \begin{multline*} #2\lnP{=}#3,\end{multline*} \else \begin{multline} \label{#1} #2\lnP{=}#3,\end{multline} \fi }
	\newcommand{\envMLineTDefPd}[3][*]{ \ifx*#1 \begin{multline*} #2\lnP{:=}#3.\end{multline*} \else \begin{multline} \label{#1} #2\lnP{:=}#3.\end{multline} \fi }
	\newcommand{\envMLineTCmDef}[3][*]{ \ifx*#1 \begin{multline*} #2\lnP{:=}#3,\end{multline*} \else \begin{multline} \label{#1} #2\lnP{:=}#3,\end{multline} \fi }
		\newcommand{\envMLinePd}{\envMLineTPd}
		\newcommand{\envMLineCm}{\envMLineTCm}

	\newcommand{\envCaseTCm}[3][?]{\begin{cases} \glcondDis{#1}#2,\lcparaCase\\\glcondDis{#1}#3,\end{cases}}
	
	\newcommand{\envCaseThCm}[4][?]{\begin{cases} \glcondDis{#1}#2,\lcparaCase\\\glcondDis{#1}#3,\lcparaCase\\\glcondDis{#1}#4,\end{cases}}
	\newcommand{\envCaseThPd}[4][?]{\begin{cases} \glcondDis{#1}#2,\lcparaCase\\\glcondDis{#1}#3,\lcparaCase\\\glcondDis{#1}#4.\end{cases}}

		\newcommand{\envLineTPdPt}[4][*]{\glcondEnvLineHead{#1} & &#3\\&#2&#4\nonumber \glcondEnvLineTailPd{#1}}

			\newcommand{\envLinePdPt}{\envLineTPdPt}

		\newcommand{\envHLineThPdPt}[5][*]{\glcondEnvLineHead{#1} #3&#2&#4\\&#2&#5\nonumber  	\glcondEnvLineTailPd{#1}}

			\newcommand{\lccondPar}[1]{\ifx#1p \\ \fi}

			\newcommand{\OTLineCThCm}[4][?]{#2&=&#3=#4,\nonumber \ifx#1p \\ \fi}
			\newcommand{\OTLineCThPd}[4][?]{#2&=&#3=#4.\nonumber \ifx#1p \\ \fi}
			\newcommand{\OTLineCThCmDef}[4][?]{#2&:=&#3=#4,\nonumber \ifx#1p \\ \fi}
			\newcommand{\OTLineCThDefPd}[4][?]{#2&:=&#3=#4.\nonumber \ifx#1p \\ \fi}
			\newcommand{\OTLineCSCm}[7][?]{#2&=&#3=#4\nonumber#5&=&#6=#7,\nonumber \ifx#1p \\ \fi}
			\newcommand{\OTLineCSPd}[7][?]{#2&=&#3=#4\nonumber#5&=&#6=#7.\nonumber \ifx#1p \\ \fi}

	\newcommand{\envMLineTCmPt}[4][*]{ \ifx*#1 \begin{multline*} #3\lnP{#2}#4,\end{multline*} \else \begin{multline} \label{#1} #3\lnP{#2}#4,\end{multline} \fi }
	\newcommand{\envMLineTPdPt}[4][*]{ \ifx*#1 \begin{multline*} #3\lnP{#2}#4.\end{multline*} \else \begin{multline} \label{#1} #3\lnP{#2}#4.\end{multline} \fi }

	\DeclareFontEncoding{OT2}{}{}
	\DeclareFontSubstitution{OT2}{cmr}{m}{n}
	\DeclareFontFamily{OT2}{cmr}{\hyphenchar\font45}
	\DeclareFontShape{OT2}{cmr}{m}{n}{<5><6><7><8><9>gen*wncyr <10><10.95><12><14.4><17.28><20.74><24.88>wncyr10}{}
	\DeclareFontShape{OT2}{cmr}{b}{n}{<5><6><7><8><9>gen*wncyb<10><10.95><12><14.4><17.28><20.74><24.88>wncyb10}{}
	\DeclareMathAlphabet{\mathcyr}{OT2}{cmr}{m}{n}
	\DeclareMathAlphabet{\mathcyb}{OT2}{cmr}{b}{n}
	\SetMathAlphabet{\mathcyr}{bold}{OT2}{cmr}{b}{n}

	\newcommand{\sh}{\mathcyr{sh}}

	\newcommand{\fcZS}[2][n]{\pwFc[#1]{\zeta}{\sh}{#2}}
	
	\newcommand{\gfcSZVs}[3][n]{\ipFc[#1]{\mathfrak{S}}{#2}{\sh}{#3}}
	\newcommand{\gfcSPL}[4][n]{\idFc[#1]{\mathfrak{SL}}{#2}{#3;#4}}

	\newcommand{\gfcDZVs}[3][n]{\ipFc[#1]{\mathfrak{D}}{#2}{\sh}{#3}}
	\newcommand{\gfcDPL}[4][n]{\idFc[#1]{\mathfrak{DL}}{#2}{#3;#4}}

	\newcommand{\gfcTZVs}[3][n]{\ipFc[#1]{\mathfrak{T}}{#2}{\sh}{#3}}
	\newcommand{\gfcTPL}[4][n]{\idFc[#1]{\mathfrak{TL}}{#2}{#3;#4}}
	\newcommand{\gfcQZV}[3][n]{\idFc[#1]{\mathfrak{Q}}{#2}{#3}}
	
	\newcommand{\gfcQZVs}[3][n]{\ipFc[#1]{\mathfrak{Q}}{#2}{\sh}{#3}}
	\newcommand{\gfcQPL}[4][n]{\idFc[#1]{\mathfrak{QL}}{#2}{#3;#4}}

	\newcommand{\lacP}[3][n]{#2\cdot#3}
	\newcommand{\lacPs}[3][n]{#2\cdot\bkS[#1]{#3}}
	\newcommand{\lacPr}[3][n]{#2\cdot\bkR[#1]{#3}}

	\newcommand{\lgpV}[1][?]{\ifx?#1V\else V_{#1}\fi}
	\newcommand{\lgpVb}[1][?]{\ifx?#1\overline{V}\else \overline{V}_{#1}\fi}
	\newcommand{\lgpD}[1][?]{\ifx?#1D\else D_{#1}\fi}

	\newcommand{\lgpVT}[2][?]{\ifx?#1{V}^{#2}\else {V}_{#1}^{#2}\fi}
	\newcommand{\lgpA}[1][?]{\ifx?#1A\else A_{#1}\fi}
	\newcommand{\lgpAb}[1][?]{\ifx?#1\overline{A}\else \overline{A}_{#1}\fi}
	\newcommand{\lgpS}[1][?]{\ifx?#1S\else S_{#1}\fi}
	\newcommand{\lgpSb}[1][?]{\ifx?#1\overline{S}\else \overline{S}_{#1}\fi}
	\newcommand{\lgpC}[1][?]{\ifx?#1C\else C_{#1}\fi}
	\newcommand{\lgpCb}[1][?]{\ifx?#1\overline{C}\else\overline{C}_{#1}\fi}
	
	\newcommand{\lsetU}[2][?]{\ifx?#1{\mathcal{U}}_{#2}\else{\mathcal{U}}^{#1}_{#2}\fi}
	\newcommand{\lsetV}[2][?]{\ifx?#1{\mathcal{V}}_{#2}\else{\mathcal{V}}^{#1}_{#2}\fi}
	\newcommand{\lsetW}[2][?]{\ifx?#1{\mathcal{W}}_{#2}\else{\mathcal{W}}^{#1}_{#2}\fi}

	\newcommand{\lmpZVi}[3][n]{\ipFc[#1]{Z}{#2}{\sh}{#3}}
	\newcommand{\lcstC}[2][b]{\Fc[#1]{C_0}{#2}}
	\newcommand{\lbfL}{{\bf l}}
	\newcommand{\lbfM}{{\bf m}}

	\newcommand{\lrpTx}[1]{(resp.\;#1)}
	
\allowdisplaybreaks[4]
	\title{\mainTitle}
	\author{\authorName}
	\date{}

\begin{document}
\maketitle

\begin{abstract}
We give a parameterized generalization of the sum formula for quadruple zeta values.
The generalization has four parameters,
	and 
	is invariant under a cyclic group of order four.
By substituting special values for the parameters,
	we also obtain weighted sum formulas for quadruple zeta values,
	which contain some known results.
\end{abstract}

\section{\sectOne} \label{sectOne}
A multiple zeta value is 
	a generalization of a classical special value of the Riemann zeta function
	\envMCm{
		\fcZ{s}
	}{
		\SmT{m=1}{\infty} \opF[s]{1}{m^s}
	}
	and is defined by 
	\envHLineDef[1_Pl_DefMZV]
	{
		\fcZ{l_1,l_2,\ldots,l_n}
	}
	{
		\Sm{m_1>m_2>\cdots>m_n>0} \opF{1}{ \pw{m_1}{l_1}\pw{m_2}{l_2}\cdots\pw{m_n}{l_n} }
	}
	for an index set $(l_1,l_2,\ldots,l_n)$ of positive integers with $l_1\geq2$.
The integers $l=l_1+\cdots+l_n$ and $n$ are called the weight and the depth respectively.
These values
	have arisen in various areas such as geometry, knot theory, number theory and mathematical physics \cite{Zagier94}.
There are many relations among these values,
	and
	an outstanding example is the sum formula
	which 
	was proved for depth two by Euler \cite{Euler1775},
	for depth three by Hoffman and Moen \cite{HM96},
	and for general depth by Granville \cite{Granville97} and Zagier \cite{Zagier95}, independently.
The formula says that the sum of all multiple zeta values of fixed weight $l$ and depth $n$ 
	is expressed by the special value $\fcZ{l}$, that is, 
	\envHLinePd[1_Pl_EqSF]
	{
		\tpSm{l_1\geq2, l_2,\ldots,l_n\geq1}{l_1+\cdots+l_n=l} \fcZ{l_1,\ldots,l_n}
	}
	{
		\fcZ{l}
	}

Various generalizations of the sum formula have been studied: 
Ohno's relations, the cyclic, restricted and weighted sum formulas
	\cite{ELO09,GX09,Ho12Ax,HO03,Machide12Ax2,Machide0000,Nakamura09,Ohno99,OZu08,OEL12Ax,SC12}.
Recently 
	parameterized generalizations of the sum formula,
	which we call parameterized sum formulas,
	were given for double and triple zeta values \cite{GKZ06,Machide12Ax1}.
The parameterized sum formula for double \lrpTx{triple} zeta values has two \lrpTx{three}  parameters,
	and is invariant under a cyclic group of order two \lrpTx{three},
	more precisely,
	the symmetry group $\gpSym{2}$ of degree two \lrpTx{the alternating group $\gpAlt{3}$ of degree three}.
By substituting special values for the parameters,
	these formulas yield some weighted sum formulas 
	which contain
	the results of Ohno and Zudilin \cite[Theorem 3]{OZu08} for double zeta values 
	and of Guo and Xie \cite[Theorem 1.1]{GX09} for triple zeta values.

In this paper, 
	we give a parameterized sum formula for quadruple zeta values
	which has four parameters
	and
	is invariant under a cyclic group of order four.
By substituting special values for the parameters,
	we also obtain weighted sum formulas for quadruple zeta values
	which contain results of Guo and Xie \cite{GX09} and of Ong, Eie and Liaw \cite{OEL12Ax}.
	
We prepare some notation in order to describe the parameterized sum formula precisely.
Let $\gpSym{n}$ be the symmetric group of degree $n$ and $\gpu$ its identity element,
	in particular, we put $\lgpS=\gpSym{4}$.
Let $\gpCyc{\sig_1, \ldots, \sig_m}$ denote the subgroup generated by permutations $\sig_1, \ldots, \sig_m$.
Let $\lgpC$ be the cyclic group $\gpCyc{(1234)}$ of order four,
	and $\lgpCb$ the subset $\bkB{\gpu, (1234)}$ of $\lgpC$,
	where $(i_1\ldots i_m)$ means a cyclic permutation defined by $i_1 \mapsto \cdots \mapsto i_m \mapsto i_1$.
Let $H_\sig$ stand for $H\sig$ for any subset $H$ and element $\sig$ of $\lgpS$.
We define a left action of $\lgpS$ on the ring $\setC[x_1,x_2,x_3,x_4]$ of polynomials in four variables
	by 
	\envMDef{
		\sig\cdot\Fc{f}{x_1,x_2,x_3,x_4}
	}{
		\Fc{f}{x_{\sig(1)},x_{\sig(2)},x_{\sig(3)},x_{\sig(4)}}
	}
	where $\sig\in\lgpS$ and $\Fc{f}{x_1,x_2,x_3,x_4}\in\setC[x_1,x_2,x_3,x_4]$.
We put
	\envM{
		x_{j_1\cdots j_m}^k
	}{
		(x_{j_1}+\cdots+x_{j_m})^k \in \setC[x_1,x_2,x_3,x_4]
	}
	for  integers $j_1,\ldots, j_m, k$ with $1\leq j_a \leq 4$, $k\geq0$. 
For example, we have for $\sig,\rho,\nu\in\lgpS$, 
	\envPLine
	{
		&&x_{1234}^k=(x_1+x_2+x_3+x_4)^k,\qquad 		\lacP{\sig\rho}{x_{234}^k}=\lacP{\sig}{x_{\rho(2)\rho(3)\rho(4)}^{k}}=x_{\sig\rho(2)\sig\rho(3)\sig\rho(4)}^{k},	\\
		&&\lacPs[b]{\sig}{ x_{1234}^{k_1} x_{234}^{k_2} +  x_{1\nu(3)}^{k_1} x_{\nu(3)2}^{k_2} }
			=
			x_{\sig(1)\sig(2)\sig(3)\sig(4)}^{k_1} x_{\sig(2)\sig(3)\sig(4)}^{k_2} +  x_{\sig(1)\sig\nu(3)}^{k_1} x_{\sig\nu(3)\sig(2)}^{k_2}.			
	}
	 
The parameterized sum formula for quadruple zeta values is as follows.
\begin{theorem}\label{1_Thm1}
Let $l$ be an integer with $l\geq5$, and $x_1,x_2,x_3,x_4$ be parameters.
Let $\pSmN[t]$ mean running over all positive integers $l_1,l_2,l_3,l_4$ satisfying $l_1\geq2$ and $l_1+l_2+l_3+l_4=l$,
	and $\fcZ{\lbfL}$ mean $\fcZ{l_1,l_2,l_3,l_4}$.
Then we have
	\envLinePd[1_Thm1_Eq1]
	{
		\pSmN
		\bkB[G]{
			\Sm{\sig\in\lgpS}  \lacPs{\sig}{x_{1234}^{l_1-1} x_{234}^{l_2-1}x_{34}^{l_3-1} x_{4}^{l_4-1}}
			\lnAHs[]{30}
			-
			\Sm{\sig\in\lgpC\cup\lgpC[(34)]}
			\lacPs[g]{\sig}{\,
				\Sm{\rho\in\gpCyc{(234)}} x_{134}^{l_1-1} x_{\rho(2)\rho(3)\rho(4)}^{l_2-1} x_{\rho(3)\rho(4)}^{l_3-1} x_{\rho(4)}^{l_4-1}
				\lnAHs{110}
				+
				\Sm{\rho\in\gpCyc{(24)}} x_{314}^{l_1-1} x_{14}^{l_2-1} x_{\rho(2)\rho(4)}^{l_3-1} x_{\rho(4)}^{l_4-1}
				+
				x_{341}^{l_1-1} x_{41}^{l_2-1} x_{1}^{l_3-1} x_{2}^{l_4-1} 
			}
			\lnAHs{30}
			+
			\Sm{\sig\in\lgpC}
			\lacPs[g]{\sig}{\,
				\Sm{\rho\in\gpCyc{(\nu(2)4)}} \Sm{\nu\in\lgpCb} x_{1\nu(3)}^{l_1-1} x_{\nu(3)2}^{l_2-1} x_{\rho\nu(2)\rho(4)}^{l_3-1} x_{\rho(4)}^{l_4-1} 
				+
				\Sm{\nu\in\lgpCb} x_{\nu(1)3}^{l_1-1} x_{2\nu(3)}^{l_2-1} x_{\nu(3)}^{l_3-1} x_{\nu(4)}^{l_4-1} 
				\lnAHs{85}
				+
				x_{41}^{l_1-1} x_{1}^{l_2-1} x_{2}^{l_3-1} x_{3}^{l_4-1} 
				-
				x_{1}^{l_1-1} x_{2}^{l_2-1} x_{3}^{l_3-1} x_{4}^{l_4-1} 
			}
		}
		\fcZ{\lbfL}
		\nonumber
	}
	{
		\bkR[G]{ \tpSm{l_1,l_2,l_3,l_4\geq1}{l_1+l_2+l_3+l_4=l} x_1^{l_1-1} x_2^{l_2-1} x_3^{l_3-1} x_4^{l_4-1} } \fcZ{l}
	}
\end{theorem}
We see that
	\refEq{1_Thm1_Eq1}  is a straightforward generalization of the original sum formula 
	because \refEq{1_Thm1_Eq1} with $(x_1,x_2,x_3,x_4)=(1,0,0,0)$ gives \refEq{1_Pl_EqSF} for $n=4$,
	and that
	\refEq{1_Thm1_Eq1} is invariant under $\lgpC$
	since
	$\Sm{\sig\in\lgpS}$ and $\Sm{\sig\in\lgpC[\alp]} (\alp\in\lgpS)$ are unchanged under any left action of $\lgpC$.

As applications of the parameterized sum formula,
	we give the weighted sum formulas in \refThm{1_Thm2} below,
	which consist of
	the known formulas \refEq{1_Thm2i_Eq1}, \refEq{1_Thm2i_Eq2} and \refEq{1_Thm2i_Eq3},
	and 
	the new formula \refEq{1_Thm2ii_Eq1} whose weights have not only powers of $2$ but also powers of $3$.
To be exact,
	\refEq{1_Thm2i_Eq1} and \refEq{1_Thm2i_Eq3} were
	the results of \cite[Theorem 1.1]{GX09} and \cite[Main Theorem]{OEL12Ax} for quadruple zeta values, respectively.
It seems that 
	\refEq{1_Thm2i_Eq2} is not printed out anywhere,
	but \refEq{1_Thm2i_Eq2} is easily derived from subtracting \refEq{1_Thm2i_Eq3} from twice \refEq{1_Thm2i_Eq1}.
\begin{theorem}[\text{cf. \cite{GX09} and \cite{OEL12Ax}}]\label{1_Thm2}
Let $l$, $\pSmN[t]$ and $\fcZ{\lbfL}$ be as in \refThm{1_Thm1}.
\mbox{}\\
{\bf(i)}
We have
	\envHLineCSNmePd
	{\label{1_Thm2i_Eq1} 
		\pSmN \bkR[b]{ 2^{l_{123}-2} + 2^{l_{12}-2} + 2^{l_1-1} - 2^{l_{23}-1} - 2^{l_2-1} } \fcZ{\lbfL}
	}
	{
		l\fcZ{l}
	}
	{\label{1_Thm2i_Eq2} 
		\pSmN \bkR[b]{ 2^{l_{123}-1} + 2^{l_{12}-1} - 2^{l_{23}} - 2^{l_2} - 2^{l_3+1} } \fcZ{\lbfL}
	}
	{
		(l-3)\fcZ{l}
	}
	{\label{1_Thm2i_Eq3}
		\pSmN \bkR[b]{ 2^{l_1} + 2^{l_3+1} }  \fcZ{\lbfL}
	}
	{
		(l+3)\fcZ{l}
	}
{\bf(ii)}
We have
	\envHLinePd[1_Thm2ii_Eq1]
	{
		\pSmN \bkR[b]{ 3^{l_2}2^{l_1-1} - 3^{l_2} - 1 } 2^{l_{13}} \fcZ{\lbfL}
	}
	{
		\opF{(l+1)(l^2+5l-18)}{12} \fcZ{l}
	}
\end{theorem}

We outline the ways to prove the theorems.
\refThm{1_Thm1} is shown by a similar way adopted in \cite{Machide12Ax1},
	that is,
	we give some identities for multiple polylogarithms instead of multiple zeta values,
	and induce \refEq{1_Thm1_Eq1} from the identities
	by using asymptotic properties of multiple polylogarithms.
Formulas \refEq{1_Thm2i_Eq1} and \refEq{1_Thm2i_Eq2} in \refThm{1_Thm2} 
	are directly proved by substituting special values for the parameters in \refEq{1_Thm1_Eq1},
	and 
	\refEq{1_Thm2i_Eq3} and \refEq{1_Thm2ii_Eq1} are induced 
	from $\setQ$-linear combinations of equations obtained by substitutions.

It is worth noting that 
	we also obtain some weighted sum formulas 
	whose weights are written in terms of powers of $2$ and $3$ in the course of the proof of \refEq{1_Thm2ii_Eq1},
	which do not have smart expressions like \refEq{1_Thm2ii_Eq1} (see \refRem{4_Rem1}). 
Furthermore,
	we prove an equation about a cyclic sum of quadruple zeta values, 
	which is related to Hoffman's result \cite[Theorem 2.2]{Hoffman92},
	in giving the identities for multiple polylogarithms (see \refRem{2.1_Rem1}).

The paper is organized as follows.
In \refSect{sectTwo} which has three subsections,
	we discuss some facts about double, triple and quadruple polylogarithms as a preparation to prove the theorems.
\refSect{sectThree} and \refSect{sectFour}
	devote the proofs of \refThm[s]{1_Thm1} and \ref{1_Thm2}, respectively.

\section{\sectTwo} \label{sectTwo}
Let $\lgP{l_1,\ldots,l_n}{z_1,\ldots,z_n}$ be the multiple polylogarithm which is defined by
	\envHLineThDef[2_Plain_DefLgP]
	{
		\lgP{l_1,\ldots,l_n}{z_1,\ldots,z_n}
	}
	{
		\Sm{m_1>\cdots>m_n>0} 
		\opF{ \pw{z_1}{m_1-m_2} \cdots \pw{z_{n-1}}{m_{n-1}-m_n}\pw{z_n}{m_n} }{ \pw{m_1}{l_1}\cdots\pw{m_{n-1}}{l_{n-1}}\pw{m_n}{l_n} }
	}
	{
		\Sm{m_1,\ldots,m_n>0} 
		\opF{ \pw{z_1}{m_1} \cdots \pw{z_{n-1}}{m_{n-1}}\pw{z_n}{m_n} }
			{ \pwR{m_1+\cdots+m_n}{l_1}\cdots\pwR{m_{n-1}+m_n}{l_{n-1}}\pw{m_n}{l_n} }
	}
	for an index set  $(l_1,\ldots,l_n)$ of positive integers and 
	a $n$-tuple $(z_1,\ldots,z_n)$ of complex numbers with $\abs{z_j}<1$.
We define parameterized sums of double, triple and quadruple polylogarithms by
	\envHLineCSCmDef[2_Plain_DefGfcMPL]
	{
		\gfcDPL{l}{x_1,x_2}{z_1,z_2}
	}
	{
		\dSmN x_1^{l_1-1}x_2^{l_2-1} \lgP{\lbfL}{z_1,z_2}
	}
	{
		\gfcTPL{l}{x_1,x_2,x_3}{z_1,z_2,z_3}
	}
	{
		\dSmN x_1^{l_1-1}x_2^{l_2-1}x_3^{l_3-1} \lgP{\lbfL}{z_1,z_2,z_3}
	}
	{
		\gfcQPL{l}{x_1,x_2,x_3,x_4}{z_1,z_2,z_3,z_4}
	}
	{
		\dSmN x_1^{l_1-1}x_2^{l_2-1}x_3^{l_3-1}x_4^{l_4-1} \lgP{\lbfL}{z_1,z_2,z_3,z_4}
	}
	where $\dSmN[t]$ means running over all positive integers $l_1,\ldots,l_n$ with $l=l_1+\cdots+l_n$
	and $\lbfL$ does $(l_1,\ldots,l_n)$ for suitable $n$.
We also denote $x^{l-1}\lgP{l}{z}$ by $\gfcSPL{l}{x}{z}$ for convenience.
We note that $\dSmN[t]$ and $\pSmN[t]$ are different;
$\dSmN[t]$ contains $l_1=1$ but $\pSmN[t]$ does not.
	
In this section, we give some identities for the above parameterized sums,
	and calculate constant terms of asymptotic expansions of functions appearing in the identities.
This section has three subsections;
The first and second subsections devote the proofs of the identities,
	which are derived from harmonic and shuffle relations (see the subsections for details of the relations).
In the third subsection, 	
	we calculate constant terms.

\subsection{\sSectTwoO} \label{sSectTwoO}
The purpose of this subsection is to prove the following identity.
We will give the proof in the end of this subsection. 
\begin{proposition}\label{2.1_Prop1}
Let $\bkA{z}^n$ denote the $n$-tuple $(z,z^2,\ldots,z^n)$,
	and $l$ be an integer with $l\geq4$. 
Then we have
	\envLinePd[2.1_Prop1_Eq]
	{
		-
		\Sm{\sig\in\lgpC} \lacP{\sig}{\gfcQPL{l}{x_1,x_2,x_3,x_4}{\bkA{z}^4}}
		\lnAH[]
		+
		\tpSm{a\geq3,b\geq1}{a+b=l} \Sm{\sig\in\lgpC} \lacPs{\sig}{ \gfcTPL{a}{x_1,x_2,x_3}{\bkA{z}^3} \gfcSPL{b}{x_4}{z} }
		\lnAH
		+
		\tpSm{a,b\geq2}{a+b=l} \Sm{\sig\in\lgpCb} \lacPs{\sig}{ \gfcDPL{a}{x_1,x_2}{\bkA{z}^2} \gfcDPL{b}{x_3,x_4}{\bkA{z}^2} }
		\lnAH
		-
		\tpSm{a\geq2,b,c\geq1}{a+b+c=l} \Sm{\sig\in\lgpC} \lacPs{\sig}{ \gfcDPL{a}{x_1,x_2}{\bkA{z}^2} \gfcSPL{b}{x_3}{z} \gfcSPL{c}{x_4}{z} }
		\lnAH
		+
		\tpSm{a,b,c,d\geq1}{a+b+c+d=l} \gfcSPL{a}{x_1}{z} \gfcSPL{b}{x_2}{z} \gfcSPL{c}{x_3}{z} \gfcSPL{d}{x_4}{z}
		\nonumber
	}
	{
		\bkR[G]{ \tpSm{l_1,l_2,l_3,l_4\geq1}{l_1+l_2+l_3+l_4=l} x_1^{l_1-1} x_2^{l_2-1} x_3^{l_3-1} x_4^{l_4-1} } \lgP{l}{z^4}
	}
\end{proposition}

We prepare some notation for the discussion after this.
Let $\lgpA\,(\subset\lgpS)$ be the alternating group of degree four.
Let $\lacPr{\sig}{f_1, \ldots, f_n}$
	denote $(\lacP{\sig}{f_1}, \ldots, \lacP{\sig}{f_n})$
	for any permutation $\sig$ of $\lgpS$ and ordered set $(f_1, \ldots, f_n)$ of polynomials in four variables.
We especially consider the case that
	each $f_j$ is expressed by $l_{j_1\cdots j_m}\ (1\leq j_a \leq 4)$,
	where $l_1,l_2,l_3,l_4$ are positive integers.
For a multiple polylogarithm $\lgP{\lbfL}{z,w,\ldots}$ having such an ordered set $\lbfL$,
	we set 
	\envMCmDef{
		\lacP{\sig}{ \lgP{\lbfL}{z,w,\ldots} }
	}{
		\lgP{\lacP{\sig}{\lbfL}}{z,w,\ldots}
	}
	and extended it to the $\setQ$-algebra spanned by $\lgP{\lbfL}{z,w,\ldots}$'s naturally.
For example,
	\envMLinePd
	{
		\lacPs{\sig}{ \lgP{l_{12},l_{34}}{z^2,z^4} + \lgP{l_{123}}{z^3}\lgP{l_4}{z} }
	\lnP{=}
		\lgP{ \lacPr{\sig}{l_{12},l_{34}} }{z^2,z^4} + \lgP{ \lacP{\sig}{l_{123}} }{z^3}\lgP{ \lacP{\sig}{l_4} }{z}
		\\
	}
	{
		\lgP{l_{\sig(1)\sig(2)},l_{\sig(3)\sig(4)}}{z^2,z^4} + \lgP{l_{\sig(1)\sig(2)\sig(3)}}{z^3}\lgP{l_{\sig(4)}}{z}
	}

The key relations for the proof of \refEq{2.1_Prop1_Eq} are the following harmonic relations 
	which are derived from decomposition of summation.
\begin{lemma}\label{2.1_Lem1}
Let $l_1,l_2,l_3,l_4$ be positive integers.
\mbox{}\\ \bfTx{(i)} 
We have 
	\envMLineCm[2.1_Lem1i_Eq]
	{
		\lgP{l_1,l_2,l_3}{\bkA{z}^3}\lgP{l_4}{z}
	}
	{
		\Sm{\sig\in\lsetU{1}} \lacP{\sig}{ \lgP{ l_1,l_2,l_3,l_4 }{\bkA{z}^4} }
		\\
		+
		\lacPs{(243)}{ \lgP{ l_{12},l_3,l_4  }{z^2,z^3,z^4} + \lgP{ l_1,l_{23},l_4 }{z,z^3,z^4} }
		+ 
		\lgP{ l_1,l_2,l_{34}  }{z,z^2,z^4} 
	}
	where 
	\envHLinePd[2.1_Lem1i_Gp]
	{ 
		 \lsetU{1}  
	}
	{  
		\bkB[n]{\gpu, (34), (243), (1432)}  
	}	
\bfTx{(ii)} 
We have
	\envMLineCm[2.1_Lem1ii_Eq]
	{
		\lgP{l_1,l_2}{\bkA{z}^2}\lgP{l_3,l_4}{\bkA{z}^2}
	}
	{
		\Sm{\sig\in\lsetU{2}} \lacP{\sig}{ \lgP{ l_1,l_2,l_3,l_4 }{\bkA{z}^4} }
		\\
		+
		\Sm{\sig\in\lsetV{2}} \lacPs{\sig}{ \lgP{ l_{12},l_3,l_4 }{z^2,z^3,z^4} + \lgP{ l_1,l_{23},l_4 }{z,z^3,z^4} + \lgP{ l_1,l_2,l_{34} }{z,z^2,z^4} }
		\\
		+
		\lacP{(23)}{ \lgP{ l_{12},l_{34} }{z^2,z^4} }
	}
	where
	\envHLinePd[2.1_Lem1ii_Gp]
	{
		\lsetV{2} 
	}
	{
		\bkB[n]{(23), (1342) },
		\quad
		\lsetU{2} 
	\lnP{=}
		\lsetV{2}  \cup \bkB[n]{\gpu, (13)(24), (132),(234)} 
	}
\bfTx{(iii)} 
We have
	\envMLineCm[2.1_Lem1iii_Eq]
	{
		\lgP{l_1,l_2}{\bkA{z}^2}\lgP{l_3}{z}\lgP{l_4}{z}
	}
	{
		\Sm{\sig\in\lsetU{3}} \lacP{\sig}{ \lgP{ l_1,l_2,l_3,l_4 }{\bkA{z}^4} }
		\\
		+
		\Sm{\sig\in\lsetV[a]{3}} \lacP{\sig}{ \lgP{ l_{12},l_3,l_4 }{z^2,z^3,z^4} }
		+
		\Sm{\sig\in\lsetV[b]{3}} \lacP{\sig}{ \lgP{ l_1,l_{23},l_4 }{z,z^3,z^4} }
		+
		\Sm{\sig\in\lsetV[c]{3}} \lacP{\sig}{ \lgP{ l_1,l_2,l_{34} }{z,z^2,z^4} }
		\\
		+
		\Sm{\sig\in\lsetW{3}} \lacP{\sig}{ \lgP{ l_{12},l_{34} }{z^2,z^4} }
		+
		\lacP{(24)}{ \lgP{ l_{123},l_4 }{z^3,z^4} }
		+
		\lgP{ l_1,l_{234} }{z,z^4} 
	}
	where
	\envHLineCSPd[2.1_Lem1iii_Gp]
	{
		\lsetW{3}
	}
	{
		\bkB[n]{(23), (24)}
		,
		\quad
		\lsetV{3}
	\lnP{=}
		\lsetW{3} \cup \bkB[n]{(34), (1342), (1423), (1432)}
	}
	{
		\lsetV[a]{3}
	}
	{
		\lsetV{3}\setminus\bkB[n]{(34)}
		,
		\quad\,
		\lsetV[b]{3}
	\lnP{=}
		\lsetV{3}\setminus\bkB[n]{(1432)}
		,
		\quad
		\lsetV[c]{3}
	\lnP{=}
		\lsetV{3}\setminus\bkB[n]{(1423)}
	}
	{
		\lsetU{3}
	}
	{
		\lsetV{3}\cup\bkB[n]{\gpu, (13)(24), (132), (142), (234), (243)}
	}
\bfTx{(iv)} 
We have
	\envMLinePd[2.1_Lem1iv_Eq]
	{
		\lgP{l_1}{z}\lgP{l_2}{z}\lgP{l_3}{z}\lgP{l_4}{z}
	}
	{
		\Sm{\sig\in\lgpS} \lacP{\sig}{ \lgP{ l_1,l_2,l_3,l_4 }{\bkA{z}^4} }
		\\
		+
		\Sm{\sig\in\lgpA} \lacPs{\sig}{ \lgP{ l_{12},l_3,l_4 }{z^2,z^3,z^4} + \lgP{ l_1,l_{23},l_4 }{z,z^3,z^4} + \lgP{ l_1,l_2,l_{34} } {z,z^2,z^4} }
		\\
		+
		\Sm{\sig\in\lgpC\cup\lgpCb[(14)]} \lacP{\sig}{ \lgP{ l_{12},l_{34} }{z^2,z^4} }
		+
		\Sm{\sig\in\lgpC} \lacPs{\sig}{ \lgP{ l_{123},l_4 }{z^3,z^4} + \lgP{ l_1,l_{234} }{z,z^4} }
		+
		\lgP{l_{1234}}{z^4}
	}
\end{lemma}
\envProof{
By \refEq{2_Plain_DefLgP},
	we have
	\envHLinePd
	{
		\lgP{l_1,l_2,l_3}{z_1,z_2,z_3}\lgP{l_4}{z_4}
	}
	{
		\Sm{m_1>m_2>m_3>0 \atop m_4>0} 
		\opF{ \pw{z_1}{m_1-m_2}\pw{z_2}{m_2-m_3}\pw{z_3}{m_3}\pw{z_4}{m_4} }{ \pw{m_1}{l_1}\pw{m_2}{l_2} \pw{m_3}{l_3} \pw{m_4}{l_4} } 
	}
From this and the decomposition of the summation 
	\envHLineCm
	{
		\Sm{m_1>m_2>m_3>0 \atop m_4>0}
	}
	{
		\Sm{m_1>m_2> m_3>m_4>0} + \Sm{m_1>m_2> m_4>m_3>0} + \Sm{m_1>m_4> m_2>m_3>0} + \Sm{m_4>m_1>m_2> m_3>0} 
		\lnAH
		+ 
		\Sm{m_1=m_4>m_2> m_3>0}  + \Sm{m_1>m_2=m_4> m_3>0} + \Sm{m_1>m_2>m_3= m_4>0}
	}
	we see that
	\envMLineCm[2.1_Lem1iP_Eq]
	{
		\lgP{l_1,l_2,l_3}{z_1,z_2,z_3}\lgP{l_4}{z_4}
	}
	{
		\lgP{l_1,l_2,l_3,l_4}{z_1,z_2,z_3,z_{34}^*} +  \lgP{l_1,l_2,l_4,l_3}{z_1,z_2,z_{24}^*,z_{34}^*} 
		\\
		+ 
		\lgP{l_1,l_4,l_2,l_3}{z_1,z_{14}^*,z_{24}^*,z_{34}^*} + \lgP{l_4,l_1,l_2,l_3}{z_4,z_{14}^*,z_{24}^*,z_{34}^*}
		\\
		+
		\lgP{l_{14},l_2,l_3}{z_{14}^*,z_{24}^*,z_{34}^*} + \lgP{l_1,l_{24},l_3}{z_1,z_{24}^*,z_{34}^*} + \lgP{l_1,l_2,l_{34}}{z_1,z_2,z_{34}^*}
	}
	where we put $z_{j_1\cdots j_m}^* = z_{j_1}\cdots z_{j_m}$.
By \refEq{2.1_Lem1iP_Eq} with $(z_1,z_2,z_3,z_4)=(z,z^2,z^3,z)$ and \refTab{Table1},
	we obtain \refEq{2.1_Lem1i_Eq}.
%
\begin{table}[!hbt]
\caption{
Each row of the lists gives 
	the conditions of $\lbfL$ and $\sig$ for $\lbfL=\lacP{\sig}{\lbfL'}$, 
	where $\lbfL' \in  \bkB{ (l_{12},l_3,l_4), (l_1,l_{23},l_4), (l_1,l_2,l_{34}), (l_{12},l_{34}), (l_{123},l_4), (l_1,l_{234})}$.
}
\label{Table1}
\begin{center}{\scriptsize\begin{tabular}{|c|c||c|c||c|c|}\hline
	\multicolumn{2}{|c||}{$\lbfL'=(l_{12},l_3,l_4)$} 	&\multicolumn{2}{|c||}{$\lbfL'=(l_1,l_{23},l_4)$} 	&\multicolumn{2}{|c|}{$\lbfL'=(l_1,l_2,l_{34})$}\\\hline
	$\lbfL$ 			     &$\sig$					&$\lbfL$ 			&$\sig$						&$\lbfL$ 			&$\sig$\\\hline	
	$(l_{12},l_3,l_4)$ &\hfill $\gpu, (12)$		&$(l_1,l_{23},l_4)$ &\hfill $\gpu, (23)$ 			&$(l_1,l_2,l_{34})$ &\hfill $\gpu, (34)$\\
	$(l_{12},l_4,l_3)$ &\hfill $(12)(34), (34)$ 	&$(l_1,l_{24},l_3)$ &\hfill $(243), (34)$ 		&$(l_1,l_3,l_{24})$ &\hfill $(234), (23)$\\
	$(l_{13},l_2,l_4)$ &\hfill $(132), (23)$ 		&$(l_1,l_{34},l_2)$ &\hfill $(234), (24)$ 		&$(l_1,l_4,l_{23})$ &\hfill $(243), (24)$\\
	$(l_{13},l_4,l_2)$ &\hfill $(234), (1342)$ 	&$(l_2,l_{13},l_4)$ &\hfill $(123), (12)$ 		&$(l_2,l_1,l_{34})$ &\hfill $(12)(34), (12)$\\
	$(l_{14},l_2,l_3)$ &\hfill $(243), (1432)$ 	&$(l_2,l_{14},l_3)$ &\hfill $(12)(34), (1243)$ 	&$(l_2,l_3,l_{14})$ &\hfill $(123), (1234)$\\
	$(l_{14},l_3,l_2)$ &\hfill $(142), (24) $ 		&$(l_2,l_{34},l_1)$ &\hfill $(124), (1234)$ 		&$(l_2,l_4,l_{13})$ &\hfill $(124), (1243)$\\
	$(l_{23},l_1,l_4)$ &\hfill $(123), (13)$ 		&$(l_3,l_{12},l_4)$ &\hfill $(132), (13)$ 		&$(l_3,l_1,l_{24})$ &\hfill $(132), (1342)$\\
	$(l_{23},l_4,l_1)$ &\hfill $(134), (1234)$ 	&$(l_3,l_{14},l_2)$ &\hfill $(13)(24), (1342)$ 	&$(l_3,l_2,l_{14})$ &\hfill $(134), (13)$\\
	$(l_{24},l_1,l_3)$ &\hfill $(143), (1243)$ 	&$(l_3,l_{24},l_1)$ &\hfill $(134), (1324)$ 		&$(l_3,l_4,l_{12})$ &\hfill $(13)(24), (1324)$\\
	$(l_{24},l_3,l_1)$ &\hfill $(124), (14)$ 		&$(l_4,l_{12},l_3)$ &\hfill $(143), (1432)$ 		&$(l_4,l_1,l_{23})$ &\hfill $(142), (1432)$\\
	$(l_{34},l_1,l_2)$ &\hfill $(13)(24), (1423)$ &$(l_4,l_{13},l_2)$ &\hfill $(142), (1423)$ 		&$(l_4,l_2,l_{13})$ &\hfill $(143), (14)$\\
	$(l_{34},l_2,l_1)$ &\hfill $(14)(23), (1324)$ &$(l_4,l_{23},l_1)$ &\hfill $(14)(23), (14)$ 		&$(l_4,l_3,l_{12})$ &\hfill $(14)(23), (1423)$\\\hline
\end{tabular}}\end{center}
\begin{center}{\scriptsize\begin{tabular}{|c|c|}\hline
	\multicolumn{2}{|c|}{$\lbfL'=(l_{12},l_{34})$} 	\\\hline
	$\lbfL$ 			     &$\sig$					\\\hline	
	$(l_{12},l_{34})$ &\hfill $\gpu, (12)(34), (12), (34)$		\\
	$(l_{13},l_{24})$ &\hfill $(132), (234), (23), (1342)$ 	\\
	$(l_{14},l_{23})$ &\hfill $(142), (243), (24), (1432)$ 		\\
	$(l_{23},l_{14})$ &\hfill $(123), (134), (13), (1234)$ 	\\
	$(l_{24},l_{13})$ &\hfill $(124), (143), (14), (1243)$ 	\\
	$(l_{34},l_{12})$ &\hfill $(13)(24), (14)(23), (1423), (1324)$ \\\hline
\end{tabular}}\end{center}
\begin{center}{\scriptsize\begin{tabular}{|c|c||c|c|}\hline
	\multicolumn{2}{|c||}{$\lbfL'=(l_{123},l_4)$} 	&\multicolumn{2}{|c|}{$\lbfL'=(l_1,l_{234})$}\\\hline
	$\lbfL$ 			     &$\sig$					&$\lbfL$ 			&$\sig$\\\hline	
	$(l_{123},l_4)$ &\hfill $\gpu, (123), (132), (12), (13), (23)$				&$(l_1,l_{234})$ &\hfill $\gpu, (234), (243), (23), (24), (34)$\\
	$(l_{124},l_3)$ &\hfill $(12)(34), (143), (243), (34), (1243), (1432)$ 	&$(l_2,l_{134})$ &\hfill $(12)(34), (123), (124), (12), (1234), (1243)$\\
	$(l_{134},l_2)$ &\hfill $(13)(24), (142), (234), (24), (1342), (1423)$ 	&$(l_3,l_{124})$ &\hfill $(13)(24), (132), (134), (13), (1324), (1342)$\\
	$(l_{234},l_1)$ &\hfill $(14)(23), (124), (134), (14), (1234), (1324)$ 	&$(l_4,l_{123})$ &\hfill $(14)(23), (142), (143), (14), (1423), (1432)$\\\hline
\end{tabular}}\end{center}
\end{table}
%

As the same way, 
	we see from the decomposition of $\Sm{m_1>m_2>0 \atop m_3>m_4>0}$ that
	\envMLineCm[2.1_Lem1iiP_Eq]
	{
		\lgP{l_1,l_2}{z_1,z_2}\lgP{l_3,l_4}{z_3,z_4}	
	}
	{
		\lgP{l_1,l_2,l_3,l_4}{z_1,z_2,z_{23}^*,z_{24}^*} + \lgP{l_1,l_3,l_2,l_4}{z_1,z_{13}^*,z_{23}^*,z_{24}^*} 
		\\
		+ 
		\lgP{l_1,l_3,l_4,l_2}{z_1,z_{13}^*,z_{14}^*,z_{24}^*} + \lgP{l_3,l_1,l_2,l_4}{z_3,z_{13}^*,z_{23}^*,z_{24}^*} +  \lgP{l_3,l_1,l_4,l_2}{z_3,z_{13}^*,z_{14}^*,z_{24}^*} 
		\\
		+ 
		\lgP{l_3,l_4,l_1,l_2}{z_3,z_4,z_{14}^*,z_{24}^*} + \lgP{l_{13},l_2,l_4}{z_{13}^*,z_{23}^*,z_{24}^*} + \lgP{l_{13},l_4,l_2}{z_{13}^*,z_{14}^*,z_{24}^*} 
		\\
		+ 
		\lgP{l_1,l_{23},l_4}{z_1,z_{23}^*,z_{24}^*} + \lgP{l_3,l_{14},l_2}{z_3,z_{14}^*,z_{24}^*} + \lgP{l_1,l_3,l_{24}}{z_1,z_{13}^*,z_{24}^*} 
		\\
		+ 
		\lgP{l_3,l_1,l_{24}}{z_3,z_{13}^*,z_{24}^*} + \lgP{l_{13},l_{24}}{z_{13}^*,z_{24}^*}
	}
	which with $(z_1,z_2,z_3,z_4)=(z,z^2,z,z^2)$ and \refTab{Table1} proves \refEq{2.1_Lem1ii_Eq}.

We will verify \refEq{2.1_Lem1iii_Eq} by using \refEq{2.1_Lem1ii_Eq}.
For this we need the following harmonic relations 
	which are obtained by the decomposition of $\Sm{m_1>0 \atop m_2>0}$ and the one of $\Sm{m_1>m_2>0 \atop m_3>0}$;
	\envHLineCFNmePd
	{\label{2.1_Lem1iiiP_Eq1}
		\lgP{k_1}{w_1}\lgP{k_2}{w_2}
	}
	{
		\lgP{k_1,k_2}{w_1,w_{12}^*} + \lgP{k_2,k_1}{w_2,w_{12}^*} + \lgP{k_{12}}{w_{12}^*}
	}
	{\label{2.1_Lem1iiiP_Eq2}
		\lgP{k_1,k_2}{w_1,w_2}\lgP{k_3}{w_3}
	}
	{
		\lgP{k_1,k_2,k_3}{w_1,w_2,w_{23}^*} 
		\lnAH[]
		+ 
		\lgP{k_1,k_3,k_2}{w_1,w_{13}^*,w_{23}^*} 
		+ 
		\lgP{k_3,k_1,k_2}{w_3,w_{13}^*,w_{23}^*} 
		\lnAH
		+ 
		\lgP{k_{13},k_2}{w_{13}^*,w_{23}^*} + \lgP{k_1,k_{23}}{w_1,w_{23}^*}
		\nonumber
	}
By \refEq{2.1_Lem1iiiP_Eq1} with $(w_1,w_2)=(z,z)$ and $(k_1,k_2)=(l_3,l_4)$,
	and \refEq{2.1_Lem1iiiP_Eq2} with $(w_1,w_2,w_3)=(z,z^2,z^2)$ and $(k_1,k_2,k_3)=(l_1,l_2,l_{34})$,
	we obtain
	\envHLineCFNmePd
	{\label{2.1_Lem1iiiP_Eq3}
		\lgP{l_3}{z}\lgP{l_4}{z}
	}
	{
		\Sm{\sig\in\gpCyc{(34)}} \lacP{\sig}{ \lgP{l_3,l_4}{\bkA{z}^2} } + \lgP{l_{34}}{z^2} 
	}
	{\label{2.1_Lem1iiiP_Eq4}
		\lgP{l_1,l_2}{\bkA{z}^2}\lgP{l_{34}}{z^2}
	}
	{
		\lacP{(1423)}{ \lgP{l_{12},l_3,l_4}{z^2,z^3,z^4} }
		\lnAH[]
		+
		\lacP{(24)}{ \lgP{l_1,l_{23},l_4}{z,z^3,z^4} } + \lacP{(34)}{ \lgP{l_1,l_2,l_{34}}{z,z^2,z^4} }
		\lnAH
		+
		\lacP{(24)}{ \lgP{l_{123},l_4}{z^3,z^4} } + \lgP{l_1,l_{234}}{z,z^4} 
		\nonumber
	}
By \refEq{2.1_Lem1iiiP_Eq3},
	we also get
	\envMLinePd[2.1_Lem1iiiP_Eq5]
	{
		\lgP{l_1,l_2}{\bkA{z}^2}\lgP{l_3}{z}\lgP{l_4}{z}
		\\
	}
	{
		\Sm{\sig\in\gpCyc{(34)}} \lacPs{\sig}{ \lgP{l_1,l_2}{\bkA{z}^2}\lgP{l_3,l_4}{\bkA{z}^2} } +  \lgP{l_1,l_2}{\bkA{z}^2}\lgP{l_{34}}{z^2}
	}
Since direct calculations show that
	\envHLineCSCm
	{
		\gpCyc{(34)}\cdot\lsetU{2}
	}
	{
		\lsetU{3}
	}
	{
		\gpCyc{(34)}\cdot\lsetV{2}
	}
	{
		\bkB[n]{(23), (1342), (243), (142) }
	}
	{
		\gpCyc{(34)}\cdot(23)
	}
	{
		\bkB[n]{(23), (243)}
	}
	\refEq{2.1_Lem1ii_Eq}, \refEq{2.1_Lem1iiiP_Eq4} and \refEq{2.1_Lem1iiiP_Eq5} with \refTab{Table1} prove \refEq{2.1_Lem1iii_Eq}.

We verify \refEq{2.1_Lem1iv_Eq} finally.
To exactly express the decomposition of a summation by permutations is difficult in general,
	but the case of $\Sm{m_1,m_2,m_3,m_4>0}$ is not
	since the summation is invariant under $\lgpS$
	and its decomposition can be understood in terms of quotient sets of $\lgpS$ as we see below.
For an odd permutation $(ij)$,
	$\lgpA$ is a transversal of $\lgpS/\gpCyc{(ij)}$
	since 
	the numbers of $\lgpA$ and $\lgpS/\gpCyc{(ij)}$ are equal
	and 
	the canonical projection of $\lgpA$ into $\lgpS/\gpCyc{(ij)}$ is injective.
It also follwos from \refTab{Table1} that
	$\lgpC\cup\lgpCb[(14)]$ is a transversal of $\lgpS/\gpCyc{(12), (34)}$,
	since
	\envMCm{
		\lgpC
	}{
		\bkB[n]{\gpu, (1234), (13)(24), (1432)}
	}
	\envMThCm{
		\lgpCb[(14)]
	}{
		\bkB[n]{\gpu, (1234) }(14)
	}{
		\bkB[n]{(14), (234)}
	}
	and
	$\lgpS/\gpCyc{(12), (34)}$
	is isomorphic to 
	$\bkB[n]{(l_{j_1j_2},l_{j_3j_4}) \mVert 1\leq j_a \leq 4, j_a\neq j_b (a\neq b) } / \! \sim$ by representation theory,
	where
	we say that
	$\lbfL_1\sim\lbfL_2$ for $\lbfL_i\in \bkB[b]{(l_{j_1j_2},l_{j_3j_4})}$ if and only if
	there is a permutation $\sig$ such that $\lbfL_1=\lacP{\sig}{\lbfL_2}$.
Similarly
	it holds that
	$\lgpC$ is a transversal of $\lgpS/\lgpS^i$, 
	where
	$\lgpS^i$ denotes $\bkB{\sig\in\lgpS \mVert \sig(i)=i}$ and $i=1,4$. 
Let $\setN$ be the set of positive integers
	and
	$\lbfM$ mean a lattice point  $(m_1,m_2,m_3,m_4)$ in $\setN^4$.
Since $\lgpA$, $\lgpC\cup\lgpCb[(14)]$ and $\lgpC$ are transversals of certain quotient sets,
	the lattice points in $\setN^4$ decompose into 
	the following disjoint subsets;
	\envPLine
	{\begin{array}{ll}
		\bkB[n]{\lbfM \mVert m_{\sig(1)}>m_{\sig(2)}>m_{\sig(3)}>m_{\sig(4)}},\quad  	& \bkB[n]{\lbfM \mVert m_{\tau(1)}=m_{\tau(2)}>m_{\tau(3)}>m_{\tau(4)}}, \\
		\bkB[n]{\lbfM \mVert m_{\tau(1)}>m_{\tau(2)}=m_{\tau(3)}>m_{\tau(4)}},\quad 	& \bkB[n]{\lbfM \mVert m_{\tau(1)}>m_{\tau(2)}>m_{\tau(3)}=m_{\tau(4)}}, \\
		\bkB[n]{\lbfM \mVert m_{\rho(1)}=m_{\rho(2)}>m_{\rho(3)}=m_{\rho(4)}},\quad	& \bkB[n]{\lbfM \mVert m_{\nu(1)}=m_{\nu(2)}=m_{\nu(3)}>m_{\nu(4)}}, \\
		\bkB[n]{\lbfM \mVert m_{\nu(1)}>m_{\nu(2)}=m_{\nu(3)}=m_{\nu(4)}},\quad		& \bkB[n]{\lbfM \mVert m_1=m_2=m_3=m_4},
	\end{array}}
	where $\sig\in\lgpS, \tau\in\lgpA, \rho\in\lgpC\cup\lgpCb[(14)]$ and $\nu\in\lgpC$. 
Similarly to \refEq{2.1_Lem1i_Eq} and \refEq{2.1_Lem1ii_Eq},
	\refEq{2.1_Lem1iv_Eq}
	follows from
	the decomposition of $\Sm{m_1,m_2,m_3,m_4>0}$ induced by the one above.
}%
We prepare some equations in order to prove \refProp{2.1_Prop1},
	which are obtained 
	by  summing up each harmonic relation in \refLem{2.1_Lem1} 
	with $(l_1,l_2,l_3,l_4)=\lacPr{\sig}{l_1,l_2,l_3,l_4}$ for all $\sig\in\lgpC$.
\begin{lemma}\label{2.1_Lem2}
Let $l_1,l_2,l_3,l_4$ be positive integers.
\mbox{}\\ \bfTx{(i)} 
We have 
	\envMLinePd[2.1_Lem2i_Eq]
	{
		\Sm{\sig\in\lgpC} \lacPs{\sig}{ \lgP{l_1,l_2,l_3}{\bkA{z}^3}\lgP{l_4}{z} }
	}
	{
		\bkR[G]{ 2\Sm{\sig\in\lgpC} + \Sm{\sig\in\lgpC[(12)]} + \Sm{\sig\in\lgpC[(34)]} } \lacP{\sig}{ \lgP{l_1,l_2,l_3,l_4}{\bkA{z}^4} }
		\\
		+
		\bkR[G]{ \Sm{\sig\in\lgpA} - \Sm{\sig\in\lgpC[(13)]} - \Sm{\sig\in\lgpC[(23)]}  }
		\lacPs{\sig}{ 
			\lgP{l_{12},l_3,l_4}{z^2,z^3,z^4} + \lgP{l_1,l_{23},l_4}{z,z^3,z^4} 
			\\
			+ 
			\lgP{l_1,l_2,l_{34}}{z,z^2,z^4} 
		}
	}
\bfTx{(ii)} 
We have
	\envMLinePd[2.1_Lem2ii_Eq]
	{
		\Sm{\sig\in\lgpCb } \lacPs{\sig}{ \lgP{l_1,l_2}{\bkA{z}^2}\lgP{l_3,l_4}{\bkA{z}^2} }
	}
	{
		\bkR[G]{ \Sm{\sig\in\lgpC} + \Sm{\sig\in\lgpC[(14)]} + \Sm{\sig\in\lgpC[(23)]} } \lacP{\sig}{ \lgP{l_1,l_2,l_3,l_4}{\bkA{z}^4} }
		\\
		+
		\Sm{\sig\in\lgpC[(23)]} \lacPs{\sig}{ \lgP{l_{12},l_3,l_4}{z^2,z^3,z^4} + \lgP{l_1,l_{23},l_4}{z,z^3,z^4} + \lgP{l_1,l_2,l_{34}}{z,z^2,z^4}  }
		\\
		+
		\Sm{\sig\in\lgpCb[(14)] } \lacP{\sig}{ \lgP{l_{12},l_{34}}{z^2,z^4} }
	}
\bfTx{(iii)} 
We have
	\envMLinePd[2.1_Lem2iii_Eq]
	{
		\Sm{\sig\in\lgpC} \lacPs{\sig}{ \lgP{l_1,l_2}{\bkA{z}^2}\lgP{l_3}{z}\lgP{l_4}{z} }
	}
	{
		\bkR[G]{ 2\Sm{\sig\in\lgpS}+\Sm{\sig\in\lgpC}-\Sm{\sig\in\lgpC[(13)]} } \lacP{\sig}{ \lgP{l_1,l_2,l_3,l_4}{\bkA{z}^4} }
		\\
		+
		\bkR[G]{ 2\Sm{\sig\in\lgpA}-\Sm{\sig\in\lgpC[(13)]} }
		\lacPs{\sig}{ \lgP{l_{12},l_3,l_4}{z^2,z^3,z^4} + \lgP{l_1,l_{23},l_4}{z,z^3,z^4} + \lgP{l_1,l_2,l_{34}}{z,z^2,z^4} }
		\\
		+
		\bkR[G]{\Sm{\sig\in\lgpC} +  2\Sm{\sig\in\lgpCb[(14)]} } \lacP{\sig}{ \lgP{l_{12},l_{34}}{z^2,z^4} }
		+
		\Sm{\sig\in\lgpC} \lacPs{\sig}{ \lgP{l_{123},l_4}{z^3,z^4} + \lgP{l_1,l_{234}}{z,z^4} }
	}
\end{lemma}
\envProof{
We see from direct calculations that
	$\lgpS$ is decomposed into the six right cosets
	\envPLine[2.1_Lem2_Eq1]
	{\begin{array}{rclrcl}
		\hspace{-15pt}	\lgpC		&=&	\bkB[n]{\gpu, (1234), (13)(24), (1432) },		\ & \lgpC[(12)]	&=&	\bkB[n]{(12), (134), (1423), (243)},  \\
		\hspace{-15pt}	\lgpC[(13)]	&=&	\bkB[n]{(13), (14)(23), (24), (12)(34)},  	 	\ & \lgpC[(14)]	&=&	\bkB[n]{(14), (234), (1243), (132)}, \\
		\hspace{-15pt}	\lgpC[(23)]	&=&	\bkB[n]{(23), (124), (1342), (143)},		 	\ & \lgpC[(34)]	&=&	\bkB[n]{(34), (123), (1324), (142)},
	\end{array}}
 	and from \refTab{Table1} that 
	\envPLine[2.1_Lem2_Eq2]
	{\begin{array}{llll}
		\lgpC\equiv\lgpC[(12)],\quad 	& \lgpC[(13)]\equiv\lgpC[(34)],\quad 	& \lgpC[(14)]\equiv\lgpC[(23)]\quad	& \text{ in } \lgpS/\gpCyc{(12)},		\\
		\lgpC\equiv\lgpC[(23)],\quad 	& \lgpC[(12)]\equiv\lgpC[(34)],\quad 	& \lgpC[(13)]\equiv\lgpC[(14)]\quad 	& \text{ in } \lgpS/\gpCyc{(23)}, 		\\
		\lgpC\equiv\lgpC[(34)],\quad 	& \lgpC[(12)]\equiv\lgpC[(13)],\quad 	& \lgpC[(14)]\equiv\lgpC[(23)]\quad	& \text{ in } \lgpS/\gpCyc{(34)}.
	\end{array}}
By \refEq{2.1_Lem2_Eq1} we obtain
	\envPLinePd
	{
		\Sm{\sig\in\lgpC\cdot\,\lsetU{1}}
	\lnP{=}
		2\Sm{\sig\in\lgpC} + \Sm{\sig\in\lgpC[(12)]} + \Sm{\sig\in\lgpC[(34)]},
		\qquad
		\Sm{\sig\in\lgpC[(243)]}
	\lnP{=}
		\Sm{\sig\in\lgpC[(12)]}
	}
Summing up \refEq{2.1_Lem1i_Eq} with $(l_1,l_2,l_3,l_4)=\lacPr{\sig}{l_1,l_2,l_3,l_4}$ for $\sig\in\lgpC$ thus gives \refEq{2.1_Lem2i_Eq}
	since $\lgpA$ is a transversal of $\lgpS/\gpCyc{(12)}, \lgpS/\gpCyc{(23)}$ or $\lgpS/\gpCyc{(34)}$.
	
Similarly we find from direct calculations that
	\envM{
		\lgpCb\cdot\lsetV{2} 
	}{
		\lgpC[(23)]
	}
	and
	\envMCm{
		\lgpCb\cdot\lsetU{2} 
	}{
		\lgpC\cup\lgpC[(14)]\cup\lgpC[(23)]
	}
	and from \refTab{Table1} that
	$\lgpCb[(23)]\equiv\lgpCb[(14)]$
	in $\lgpS/\gpCyc{(12), (34)}$.
Therefore summing up \refEq{2.1_Lem1ii_Eq} with $(l_1,l_2,l_3,l_4)=\lacPr{\sig}{l_1,l_2,l_3,l_4}$ for $\sig\in\lgpCb$ gives \refEq{2.1_Lem2ii_Eq}.

We see from \refEq{2.1_Lem2_Eq1} that
	\envHLineCSCm
	{
		\Sm{\sig\in\lgpC\cdot\lsetW{3}}
	}
	{
		\Sm{\sig\in\lgpC[(13)]} + \Sm{\sig\in\lgpC[(23)]}
	}
	{
		\Sm{\sig\in\lgpC\cdot\lsetV{3}}
	}
	{
		\Sm{\sig\in\lgpC} + \Sm{\sig\in\lgpC[(12)]} + \Sm{\sig\in\lgpC[(13)]} + 2\Sm{\sig\in\lgpC[(23)]} + \Sm{\sig\in\lgpC[(34)]} 
	}
	{
		\Sm{\sig\in\lgpC\cdot\,\lsetU{3}}
	}
	{
		2\Sm{\sig\in\lgpS}+\Sm{\sig\in\lgpC}-\Sm{\sig\in\lgpC[(13)]}
	}
	and from \refTab{Table1} that 
	$\lgpC[(13)]\equiv\lgpC$, $\lgpC[(23)]\equiv\lgpCb[(14)]$ in $\lgpS/\gpCyc{(12), (34)}$
	and
	$\lgpC[(24)]\equiv\lgpC$ in $\lgpS/\lgpS^4$.
Because of \refEq{2.1_Lem2_Eq2},
	summing up \refEq{2.1_Lem1iii_Eq} with $(l_1,l_2,l_3,l_4)=\lacPr{\sig}{l_1,l_2,l_3,l_4}$ for $\sig\in\lgpC$ also gives \refEq{2.1_Lem2iii_Eq}.
}
We prove \refProp{2.1_Prop1}.
\envProof[\refProp{2.1_Prop1}]{
It follows from \refEq{2.1_Lem1iv_Eq}, \refEq{2.1_Lem2i_Eq}, \refEq{2.1_Lem2ii_Eq} and \refEq{2.1_Lem2iii_Eq} that
	\envLinePd[2.1_Prop1P_Eq1]
	{
		\Sm{\sig\in\lgpC} \lacPs{\sig}{ \lgP{l_1,l_2,l_3}{\bkA{z}^3}\lgP{l_4}{z} - \lgP{l_1,l_2}{\bkA{z}^2}\lgP{l_3}{z}\lgP{l_4}{z} }
		\lnAH
		+
		\Sm{\sig\in\lgpCb } \lacPs{\sig}{ \lgP{l_1,l_2}{\bkA{z}^2}\lgP{l_3,l_4}{\bkA{z}^2} }
		+
		\lgP{l_1}{z}\lgP{l_2}{z}\lgP{l_3}{z}\lgP{l_4}{z}
	}
	{
		\Sm{\sig\in\lgpC}  \lacP{\sig}{ \lgP{l_1,l_2,l_3,l_4}{\bkA{z}^4} } + \lgP{l_{1234}}{z^4}
	}
Adding together \refEq{2.1_Prop1P_Eq1} up to $x_1^{l_1-1}x_2^{l_2-1}x_3^{l_3-1}x_4^{l_4-1}$ 
	for all positive integers $l_1,l_2,l_3,l_4$ with $l_1+l_2+l_3+l_4=l$
	proves \refEq{2.1_Prop1_Eq}.
}
\begin{remark}\label{2.1_Rem1}
Let $\lbfL=(l_1, l_2, l_3, l_4)$ be an index set of integers with $l_j\geq2$.
By \refEq{2.1_Prop1P_Eq1} with $z=1$,
	we obtain an equation for a cyclic sum of $\fcZ{\lbfL}$,
	\envMLinePd[2.1_Rem1_EqMine]
	{
		\Sm{\sig\in\lgpC} \lacP{\sig}{ \fcZ{\lbfL} }
	}
	{
		\fcZ{l_1}\fcZ{l_2}\fcZ{l_3}\fcZ{l_4}
		+
		\Sm{\sig\in\lgpCb} \lacPs{\sig}{ \fcZ{l_1, l_2}\fcZ{l_3, l_4} }
		\\
		+
		\Sm{\sig\in\lgpC} \lacPs{\sig}{ \fcZ{l_1, l_2, l_3}\fcZ{l_4} - \fcZ{l_1,l_2}\fcZ{l_3}\fcZ{l_4} }
		-
		\fcZ{l_{1234}}
	}
On the other hand,
	Hoffman \cite[Theorem 2.2]{Hoffman92} proved equations for symmetric sums of multiple zeta values.
The equations in cases of double, triple and quadruple zeta values are as follows.
(Note that $A$ instead of $\fcZ[]{}$ is used for a symbol of multiple zeta value in \cite{Hoffman92}.)
	\envHLineCSCmNme
	{\label{2.1_Rem1_Eq1Hoffman}
		\Sm{\sig\in\gpSym{2}} \lacP{\sig}{\fcZ{\lbfL_2}}
	}
	{
		\fcZ{l_1}\fcZ{l_2}-\fcZ{l_{12}}
	}
	{\label{2.1_Rem1_Eq2Hoffman}
		\Sm{\sig\in\gpSym{3}} \lacP{\sig}{\fcZ{\lbfL_3}}
	}
	{
		\fcZ{l_1}\fcZ{l_2}\fcZ{l_3} - \Sm{\sig\in\gpCyc{(123)}} \lacPs{\sig}{ \fcZ{l_{12}}\fcZ{l_3} } + 2\fcZ{l_{123}}
	}
	{\label{2.1_Rem1_Eq3Hoffman}
		\Sm{\sig\in\lgpS} \lacP{\sig}{ \fcZ{\lbfL} }
	}
	{
		\fcZ{l_1}\fcZ{l_2}\fcZ{l_3}\fcZ{l_4}
		-
		\Sm{\sig\in\lgpC\cup\lgpCb[(14)]} \lacPs{\sig}{ \fcZ{l_{12}}\fcZ{l_3}\fcZ{l_4} }
		\lnAHs[]{-13}
		+
		\Sm{\sig\in\gpCyc{(123)}} \lacPs{\sig}{ \fcZ{l_{12}}\fcZ{l_{34}} }
		+
		2\Sm{\sig\in\lgpC} \lacPs{\sig}{ \fcZ{l_{123}}\fcZ{l_4} }
		-
		6\fcZ{l_{1234}}
		\nonumber
	}
	where we put $\lbfL_2=(l_1,l_2)$ and $\lbfL_3=(l_1,l_2,l_3)$.

By virtue of \refEq{2.1_Rem1_Eq1Hoffman} and \refEq{2.1_Rem1_Eq2Hoffman},
	summing up  \refEq{2.1_Rem1_EqMine} with the left actions of $\sig\in\bkB{\gpu, (12), (13), (14), (23), (34)}$
	yields \refEq{2.1_Rem1_Eq3Hoffman} after some calculations.
\end{remark}

\subsection{\sSectTwoT} \label{sSectTwoT}
The purpose of this subsection is to prove the following identities.
\begin{proposition}\label{2.2_Prop1}
Let $\pwB{z}{n}$ denote the $n$-tuple $(z,\ldots,z)$,
	and $l$ be an integer with $l\geq4$.
\mbox{}\\ \bfTx{(i)} 
We have
	\envMLinePd[2.2_Prop1i_Eq]
	{
		\tpSm{a\geq3,b\geq1}{a+b=l} \gfcTPL{a}{x_{1},x_{2},x_{3}}{\pwB{z}{3}} \gfcSPL{b}{x_4}{z}
	}
	{
		\Sm{\rho\in\gpCyc{(34)}} \gfcQPL{l}{x_{14},x_{24},x_{\rho(3)\rho(4)},x_{\rho(4)}}{\pwB{z}{4}}
		\\
		+
		\gfcQPL{l}{x_{14},x_{42},x_{2},x_{3}}{\pwB{z}{4}}
		+
		\gfcQPL{l}{x_{41},x_{1},x_{2},x_{3}}{\pwB{z}{4}}
	}
\bfTx{(ii)} 
We have
	\envLinePd[2.2_Prop1ii_Eq]
	{
		\tpSm{a,b\geq2}{a+b=l} \gfcDPL{a}{x_{1},x_{2}}{\pwB{z}{2}} \gfcDPL{b}{x_{3},x_{4}}{\pwB{z}{2}}
	}
	{
		\Sm{\sig\in\gpCyc{(13)(24)}} 
		\lacPs[g]{\sig}{ 
			\Sm{\rho\in\gpCyc{(24)}}  \gfcQPL{l}{x_{13},x_{32},x_{\rho(2)\rho(4)},x_{\rho(4)}}{\pwB{z}{4}} 
			+ 
			\gfcQPL{l}{x_{13},x_{23},x_{3},x_{4}}{\pwB{z}{4}}
		}
	}
\bfTx{(iii)} 
We have
	\envLinePd[2.2_Prop1iii_Eq]
	{
		\tpSm{a\geq2,b,c\geq1}{a+b+c=l} \gfcDPL{a}{x_{1},x_{2}}{\pwB{z}{2}} \gfcSPL{b}{x_3}{z} \gfcSPL{c}{x_4}{z}
	}
	{
		\Sm{\sig\in\gpCyc{(34)}}
		\lacPs[g]{\sig}{
			\Sm{\rho\in\gpCyc{(234)}} \gfcQPL{l}{x_{134},x_{\rho(2)\rho(3)\rho(4)},x_{\rho(3)\rho(4)},x_{\rho(4)}}{\pwB{z}{4}}
			\lnAHs{30}
			+
			\Sm{\rho\in\gpCyc{(24)}} \gfcQPL{l}{x_{314},x_{14},x_{\rho(2)\rho(4)},x_{\rho(4)}}{\pwB{z}{4}}
			+
			\gfcQPL{l}{x_{341},x_{41},x_{1},x_{2}}{\pwB{z}{4}}
		}
	}
\bfTx{(iv)} 
We have
	\envLinePd[2.2_Prop1iv_Eq]
	{
		\tpSm{a,b,c,d\geq1}{a+b+c+d=l} \gfcSPL{a}{x_1}{z} \gfcSPL{b}{x_2}{z} \gfcSPL{c}{x_3}{z} \gfcSPL{d}{x_4}{z}
	}
	{
		\Sm{\sig\in\lgpS} \lacP{\sig}{ \gfcQPL{l}{x_{1234},x_{234},x_{34},x_{4}}{\pwB{z}{4}} }
	}
\end{proposition}
As we see in the proof below,
	each coefficient of $x_1^{l_1-1}x_2^{l_2-1}x_3^{l_3-1}x_4^{l_4-1}$'s of the above identities
	is a shuffle relation for quadruple polylogarithms 
	because of the equivalence between partial fraction expansions and shuffle relations (see \cite{GX11} and \cite{KMT11}).
\envProof[\refProp{2.2_Prop1}]{
By the partial fraction procedure (see \cite[Lemma 2.3]{Yamamoto11Ax} for example), we obtain
	\envHLine
	{
		\opF{1}{m_{123}m_{23}m_3 m_4} 
	}
	{
		\Sm{\sig\in \lsetU{1} } \opF{1}{ m_{\sig(1)\sig(2)\sig(3)\sig(4)}m_{\sig(2)\sig(3)\sig(4)}m_{\sig(3)\sig(4)}m_{\sig(4)} }
	}
	where $\lsetU{1}$ is the set defined in \refEq{2.1_Lem1i_Gp}.
By replacing $m_j$ by $m_j-ty_j$
	and calculating the coefficient of $t^{l-4}$,
	we get
	\envMLinePd[2.2_Prop1PEq1]
	{
		\tpSm{a\geq3, b\geq1}{a+b=l}
		\bkR[G]{
			\tpSm{l_1,l_2,l_3\geq1}{l_1+l_2+l_3=a} \opF{ y_{123}^{l_1-1}y_{23}^{l_2-1}y_{3}^{l_3-1} }{ m_{123}^{l_1}m_{23}^{l_2}m_{3}^{l_3} }
		}
		\opF{ y_4^{b-1} }{ m_4^{b} }
		\\
	}
	{
		\Sm{\sig\in \lsetU{1} } \tpSm{l_1,l_2,l_3,l_4\geq1}{l_1+l_2+l_3+l_4=l} 
		\opF{ y_{\sig(1)\sig(2)\sig(3)\sig(4)}^{l_1-1}y_{\sig(2)\sig(3)\sig(4)}^{l_2-1}y_{\sig(3)\sig(4)}^{l_3-1}y_{\sig(4)}^{l_4-1} }
			{ m_{\sig(1)\sig(2)\sig(3)\sig(4)}^{l_1}m_{\sig(2)\sig(3)\sig(4)}^{l_2}m_{\sig(3)\sig(4)}^{l_3}m_{\sig(4)}^{l_4} }
	}
By \refEq{2_Plain_DefLgP} and \refEq{2_Plain_DefGfcMPL}, 
the sum of \refEq{2.2_Prop1PEq1} up to $z^{m_{1234}}$ for all positive integers $m_1, m_2, m_3, m_4$
	gives
	\envHLineCm
	{
		\tpSm{a\geq3,b\geq1}{a+b=l} \gfcTPL{a}{y_{123},y_{23},y_{3}}{\pwB{z}{3}} \gfcSPL{b}{y_4}{z}
	}
	{
		\Sm{\sig\in \lsetU{1} } \lacP{\sig}{ \gfcQPL{l}{y_{1234},y_{234},y_{34},y_{4}}{\pwB{z}{4}} }
	}
	which with $(y_1,y_2,y_3,y_4)=(x_1-x_2,x_2-x_3,x_3,x_4)$ proves \refEq{2.2_Prop1i_Eq}.

In the same way,
	we can find from the partial fraction expansion
	\envHLine
	{
		\opF{1}{m_{12}m_2 m_{34}m_4 }
	}
	{
		\Sm{\sig\in\lsetU{2} } \opF{1}{ m_{\sig(1)\sig(2)\sig(3)\sig(4)}m_{\sig(2)\sig(3)\sig(4)}m_{\sig(3)\sig(4)}m_{\sig(4)} }
	}
	that
	\envHLineCm
	{
		\tpSm{a,b\geq2}{a+b=l} \gfcDPL{a}{y_{12},y_2}{\pwB{z}{2}} \gfcDPL{b}{y_{34},y_4}{\pwB{z}{2}}
	}
	{
		\Sm{\sig\in \lsetU{2} } \lacP{\sig}{ \gfcQPL{l}{y_{1234},y_{234},y_{34},y_{4}}{\pwB{z}{4}} }
	}
	which with $(y_1,y_2,y_3,y_4)=(x_1-x_2,x_2,x_3-x_4,x_4)$ proves \refEq{2.2_Prop1ii_Eq}.

Identity \refEq{2.2_Prop1iii_Eq} can be reduced to \refEq{2.2_Prop1ii_Eq} by using
	\envHLine
	{
		\tpSm{b,c\geq1}{b+c=k} \gfcSPL{b}{x_3}{z} \gfcSPL{c}{x_4}{z}
	}
	{
		\Sm{\sig\in\gpCyc{(34)}} \lacP{\sig}{ \gfcDPL{k}{x_{34},x_4}{\pwB{z}{2}} }
	}
	for integers $k\geq2$,
	which are obtained by the typical partial fraction expansion 
	$\opF[s]{1}{m_3m_4}=\opF[s]{1}{m_{34}m_4} + \opF[s]{1}{m_{43}m_3}$. 

Identity \refEq{2.2_Prop1iv_Eq} is proved by the partial fraction expansion
	\envHLine
	{
		\opF{1}{m_1m_2m_3m_4}
	}
	{
		\Sm{\sig\in \lgpS  } \opF{1}{ m_{\sig(1)\sig(2)\sig(3)\sig(4)}m_{\sig(2)\sig(3)\sig(4)}m_{\sig(3)\sig(4)}m_{\sig(4)} }
	}
	similarly to \refEq{2.2_Prop1i_Eq} and \refEq{2.2_Prop1ii_Eq}.
}

\subsection{\sSectTwoTh} \label{sSectTwoTh}
Let $\sbLandau[]{}$ denote the Landau symbol.
For any function $\Fc{F}{z}$ which has a polynomial
	$\Fc{P}{T}$ and a positive number $J>0$, and satisfies the asymptotic property 
	\envHLineCm[2.3_Plain_EqAsymptoticForm]
	{
		\Fc{F}{z}
	}
	{
		\Fc[b]{P}{-\lgg[n]{1-z}}
		+
		\sbLandau[B]{(1-z)(\lgg[n]{1-z})^J}
		\quad
		(z\nearrow1)
	}
	we define the constant term of $\Fc{P}{T}$ by $\lcstC{\Fc{F}{z}}$,
	that is,
	\envMPd
	{
		\lcstC{\Fc{F}{z}}
	}
	{
		\Fc{P}{0}
	}
Here $z\nearrow1$ means bringing $z$ close to $1$ under the condition $0<z<1$.
Function $\lcstC[]{}$ is well defined
	since $\Fc{P}{T}$ is uniquely determined by \refEq{2.3_Plain_EqAsymptoticForm}.
	
In this subsection, 
	we evaluate 
	images of some functions written in terms of multiple polylogarithms under $\lcstC[]{}$,
	which enable us to calculate the constant terms of the asymptotic expansions of the functions 
	appearing in \refProp[s]{2.1_Prop1} and \ref{2.2_Prop1}.

Firstly we introduce the basic asymptotic properties about the multiple polylogarithms $\lgP{l_1,\ldots,l_n}{\pwB{z}{n}}$ 
	which were shown in \cite[\S2]{IKZ06}.
\begin{lemma}[\cite{IKZ06}]\label{2.3_LemmaA}
For any multiple polylogarithm $\lgP{l_1,\ldots,l_n}{\pwB{z}{n}}$,
	there are 
	a unique polynomial $\lmpZVi{l_1,\ldots,l_n}{T}$ and a positive number $J>0$
	such that \refEq{2.3_Plain_EqAsymptoticForm} holds.
\end{lemma}
We call $\lmpZVi{l_1,\ldots,l_n}{0}$ the regularized multiple zeta value, 
	and denote it by $\fcZS{l_1,\ldots,l_n}$.
It is obvious by definition that
	\envHLinePd[2.3_Plain_EqCofMPL]
	{
		\lcstC{ \lgP{l_1,\ldots,l_n}{\pwB{z}{n}} }
	}
	{
		\fcZS{l_1,\ldots,l_n}
	}
For example, 
	$\fcZS{l_1,\ldots,l_n}=\fcZ{l_1,\ldots,l_n}$ for $l_1\geq2$,
	and
	\envM
	{
		\fcZS{ \bkB{1}^n }
	}
	{
		0
	}
	because
	\envPLineCm[2.3_Plain_EqExOfSHRel]
	{
		\lgP{1}{z}
	\lnP{=}
		-\lgg[n]{1-z}
		\quad\text{ and }\quad
		\lgP{\pwB{1}{n}}{\pwB{z}{n}}
	\lnP{=}
		\opF{ \lgP{1}{z}^n }{n!}
	}
	where the last equation is a shuffle relation
	derived from
	\envHLinePd
	{
		\opF{1}{m_1\cdots m_n}
	}
	{
		\Sm{\sig\in\gpSym{n}} \opF{1}{m_{\sig(1)\sig(2)\cdots\sig(n)} m_{\sig(2)\cdots\sig(n)} \cdots m_{\sig(n)}}
	}

We define parameterized sums of regularized double, triple and quadruple zeta values by
	\envHLineCSCmDef[2.3_Plain_DefGfcRZV]
	{
		\gfcDZVs{l}{x_1,x_2}
	}
	{
		\dSmN x_1^{l_1-1}x_2^{l_2-1} \fcZS{\lbfL}
	}
	{
		\gfcTZVs{l}{x_1,x_2,x_3}
	}
	{
		\dSmN x_1^{l_1-1}x_2^{l_2-1}x_3^{l_3-1} \fcZS{\lbfL}
	}
	{
		\gfcQZVs{l}{x_1,x_2,x_3,x_4}
	}
	{
		\dSmN x_1^{l_1-1}x_2^{l_2-1}x_3^{l_3-1}x_4^{l_4-1} \fcZS{\lbfL}
	}
	where $\lbfL$ means $(l_1,\ldots,l_n)$ for suitable $n$.
We also denote $x^{l-1}\fcZS{l}$ by $\gfcSZVs{l}{x}$ for convenience,
	and easily see that
	\envM{
		\lcstC{ \gfcSPL{l}{x}{z} }
	}{
		\gfcSZVs{l}{x}
	}
	by definition.
	
The objective images under $\lcstC[]{}$ are in \refProp[s]{2.3_Prop1} and \ref{2.3_Prop2} below.
The images in \refProp{2.3_Prop1} \lrpTx{\ref{2.3_Prop2}}
	enable us 
	to calculate the constant terms of the asymptotic expansions of the functions appearing in \refProp{2.1_Prop1} \lrpTx{\ref{2.2_Prop1}}.
\begin{proposition}\label{2.3_Prop1}
Let $l$ be a positive integer.
We assume that
	$l\geq2$ in the first,
	$l\geq3$ in the second,
	and
	$l\geq5$ in the third equation.
Then we have
	\envHLineCSNme
	{\label{2.3_Prop1_Eq1}
		\lcstC{ \gfcDPL{l}{x_1,x_2}{\bkA{z}^2} }
	}
	{
		\gfcDZVs{l}{x_1,x_2}
		+
		\envCaseTCm[d]{ 
			0 &(l>2)
		}{
			- \opF{\fcZ{2}}{2} & (l=2)
		}
	}
	{\label{2.3_Prop1_Eq2}
		\lcstC{ \gfcTPL{l}{x_1,x_2,x_3}{\bkA{z}^3} }
	}
	{
		\gfcTZVs{l}{x_1,x_2,x_3}
		\lnAH
		+
		\envCaseTCm[d]{
			- \opF{\fcZ{2}\fcZ{l-2}}{2} x_3^{l-3} & (l>3)
		}{
			\opF{\fcZ{3}}{3} & (l=3)
		}
	}
	{\label{2.3_Prop1_Eq3}
		\lcstC{ \gfcQPL{l}{x_1,x_2,x_3,x_4}{\bkA{z}^4} }
	}
	{
		\gfcQZVs{l}{x_1,x_2,x_3,x_4} 
		\lnAH
		- 
		\opF{\fcZ{2}}{2} \gfcDZVs{l-2}{x_3,x_4} 
		+ 
		\opF{\fcZ{3}\fcZ{l-3}}{3} x_4^{l-4}
		.
	}
\end{proposition}
\begin{proposition}\label{2.3_Prop2}
Let $l$ be an integer as in \refProp{2.3_Prop1}.
Then we have
	\envHLineCSNmePd
	{\label{2.3_Prop2_Eq1}
		\lcstC{ \gfcDPL{l}{x_1,x_2}{\pwB{z}{2}} }
	}
	{
		\gfcDZVs{l}{x_1,x_2}
	}
	{\label{2.3_Prop2_Eq2}
		\lcstC{ \gfcTPL{l}{x_1,x_2,x_3}{\pwB{z}{3}} }
	}
	{
		\gfcTZVs{l}{x_1,x_2,x_3}
	}
	{\label{2.3_Prop2_Eq3}
		\lcstC{ \gfcQPL{l}{x_1,x_2,x_3,x_4}{\pwB{z}{4}} }
	}
	{
		\gfcQZVs{l}{x_1,x_2,x_3,x_4}
	}	
\end{proposition}

By comparing \refProp[s]{2.3_Prop1} with \ref{2.3_Prop2},
	we see that the images of 
	\envPLine{
		\text{$\gfcDPL{l}{x_1,x_2}{\bkA{z}^2}$, $\gfcTPL{l}{x_1,x_2,x_3}{\bkA{z}^3}$ and $\gfcQPL{l}{x_1,x_2,x_3,x_4}{\bkA{z}^4}$}
	}
	under $\lcstC[]{}$
	are nearly equal to those of 
	\envPLineCm{
		\text{$\gfcDPL{l}{x_1,x_2}{\pwB{z}{2}}$, $\gfcTPL{l}{x_1,x_2,x_3}{\pwB{z}{3}}$ and $\gfcQPL{l}{x_1,x_2,x_3,x_4}{\pwB{z}{4}}$}
	}
	respectively,
	but they have some extra values.

We will prove the two propositions.
\refProp{2.3_Prop2} easily follows from \refEq{2_Plain_DefGfcMPL}, \refEq{2.3_Plain_EqCofMPL} and \refEq{2.3_Plain_DefGfcRZV}.
We give two lemmas
	in order to prove \refProp{2.3_Prop1}.
	
\begin{lemma}\label{2.3_Lem1}
Let $(l_1,\ldots,l_n), (i_1,\ldots,i_n)$ and $(j_1,\ldots,j_n)$
	be $n$-tuples of positive integers. 
Assume that
	$i_1\leq\cdots\leq i_n$, $j_1\leq\cdots\leq j_n$ and $i_a\leq j_a$ for every integer $a$. 
If there is a positive integer $d$ such that $d<n$, $l_1=\cdots=l_d=1$, $l_{d+1}>1$ and $i_a=j_a\ (a=1,\ldots,d)$,
	then
	\envHLinePd[2.3_Lem1_Eq]
	{
		\Lm{z\nearrow1}
		\bkR{ \lgP{l_1,\ldots,l_n}{z^{i_1}, \ldots, z^{i_n}} - \lgP{l_1,\ldots,l_n}{z^{j_1}, \ldots, z^{j_n}} }
	}
	{
		0
	}
\end{lemma}
\envProof{
If $i_a=j_a$ for any integer $a$ with $d+1\leq a \leq n$, 
	then $\lgP{l_1,\ldots,l_n}{z^{i_1}, \ldots, z^{i_n}} = \lgP{l_1,\ldots,l_n}{z^{j_1}, \ldots, z^{j_n}}$
	and
	\refEq{2.3_Lem1_Eq} holds evidently.
Suppose that
	it is false
	and $b$ is the smallest integer such that $i_b<j_b$.
We set $i_0=j_0=0$,
	and 
	put
	$i_a'=i_a-i_{a-1}$ and $j_a'=j_a-j_{a-1}$ for every integer $a$.
From the conditions $0<i_1\leq\cdots\leq i_n$ and $0<j_1\leq\cdots\leq j_n$,
	we easily see that $i_1', j_1'>0$ and $i_a', j_a'\geq0$ for any integer $a$ with $2\leq a \leq n$.
Since
	$i_a'=j_a'\ (a=1,\ldots,b-1)$ and $i_b'<j_b'$ by the minimality of $b$,
	we have
	\envLineThPd
	{
		\lgP{l_1,\ldots,l_n}{z^{i_1}, \ldots, z^{i_n}} - \lgP{l_1,\ldots,l_n}{z^{j_1}, \ldots, z^{j_n}}
	}
	{
		\Sm{m_1>\cdots>m_n>0} \opF{ z^{i_1'm_1+\cdots+i_n'm_n} - z^{j_1'm_1+\cdots+j_n'm_n} }{m_1^{l_1}\cdots m_n^{l_n}}
	}
	{
		\Sm{m_1>\cdots>m_n>0} 
		\opF{ z^{i_1'm_1+\cdots+i_b'm_b} (z^{i_{b+1}'m_{b+1}+ \cdots + i_n'm_n} -z^{ (j_b'-i_b')m_b+ j_{b+1}'m_{b+1}+ \cdots + j_n'm_n })  }
			{m_1^{l_1}\cdots m_n^{l_n}}
	}
We put $i=i_1'$ and $j=\Max{j_b'-i_b', j_{b+1}', \ldots, j_n'}$.
It is clear that $i,j>0$.
Under the condition $0<z<1$,
	we thus see that
	\envHLineThPdPt[2.3_Lem1P_Eq1]{<}
	{
		0
	}
	{
		\lgP{l_1,\ldots,l_n}{z^{i_1}, \ldots, z^{i_n}} - \lgP{l_1,\ldots,l_n}{z^{j_1}, \ldots, z^{j_n}}
	\lnAHP[]{\leq}
		\Sm{m_1>\cdots>m_n>0}  \opF{ z^{im_1} (1 -z^{ (j_b'-i_b')m_b+ j_{b+1}'m_{b+1}+ \cdots + j_n'm_n })  }{m_1^{l_1}\cdots m_n^{l_n}}
	\lnAHP{\leq}
		\Sm{m_1>\cdots>m_n>0}  \opF{ z^{im_1} (1 -z^{ (j_b'-i_b'+j_{b+1}'+\cdots+j_n')m_{d+1} })  }{m_1^{l_1}\cdots m_n^{l_n}}
		\nonumber
	}
	{
		\Sm{m_1>\cdots>m_n>0} 
		\opF{ z^{im_1} (1 -z^{ njm_{d+1}}) }{m_1^{l_1}\cdots m_n^{l_n}}
	}
Since
	$1-z^k=(1-z)\SmT{h=0}{k-1}z^h<k(1-z)$
	and
	$l_{d+1}>1$,
	it follows that
	\envLinePdPt[2.3_Lem1P_Eq2]{<}
	{
		\Sm{m_1>\cdots>m_n>0} \opF{ z^{im_1} (1 -z^{ njm_{d+1} })  }{m_1^{l_1}\cdots m_n^{l_n}}
	}
	{
		nj(1-z) \Sm{m_1>\cdots>m_n>0} \opF{ z^{im_1} }{m_1^{l_1}\cdots m_{d+1}^{l_{d+1}-1}\cdots m_n^{l_n}}
	\lnAHP{\leq}
		nj(1-z) \Sm{m_1>\cdots>m_n>0} \opF{ z^{im_1} }{m_1\cdots m_n}
	\lnAHP{=}
		nj(1-z) \lgP{\pwB{1}{n}}{ \bkB[n]{z^{i}}^n } 
	}
We see from \refEq{2.3_Plain_EqExOfSHRel} that
	\envHLineCm
	{
		\Lm{z\nearrow1} (1-z) \lgP{\pwB{1}{n}}{ \bkB[n]{z^{i}}^n } 
	}
	{
		0
	}
	which with \refEq{2.3_Lem1P_Eq1} and \refEq{2.3_Lem1P_Eq2} proves \refEq{2.3_Lem1_Eq}.
}
We note the following facts;
When $l_1,\ldots,l_n$ are positive integers,
	$l_1\cdots l_n=1$ if and only if $l_1=\cdots =l_n=1$,
	and 
	$l_1\cdots l_n>1$ if and only if $l_j>1$ for some $j$. 
\begin{lemma}\label{2.3_Lem2}
Let $l_1,l_2,l_3,l_4,l$ be positive integers.
We assume that
	$l=l_1+l_2+l_3$ in the second,
	and
	$l=l_1+l_2+l_3+l_4$ and $l\geq5$ in the third equation.
Then we have
	\envHLineCSNme
	{\label{2.3_Lem2_Eq1}
		\lcstC{ \lgP{l_1,l_2}{\bkA{z}^2} }
	}
	{
		\fcZS{l_1,l_2}
		+
		\envCaseTCm{
			0 & (l_1l_2>1)
		}{
			- \opF{\fcZ{2}}{2} & (l_1l_2=1)
		}
	}
	{\label{2.3_Lem2_Eq2}
		\lcstC{ \lgP{l_1,l_2,l_3}{\bkA{z}^3} }
	}
	{
		\fcZS{l_1,l_2,l_3}
		\lnAH
		+
		\envCaseThCm{
			0 & (l_1l_2>1)
		}{
			- \opF{\fcZ{2}\fcZ{l-2}}{2} & (l_1l_2=1 \text{\;and\ } l_3>1)
		}{
			\opF{\fcZ{3}}{3}  & (l_1l_2l_3=1)
		}
	}
	{\label{2.3_Lem2_Eq3}
		\lcstC{ \lgP{l_1,l_2,l_3,l_4}{\bkA{z}^4} }
	}
	{
		\fcZS{l_1,l_2,l_3,l_4}
		\lnAH
		+
		\envCaseThPd{
			0 & (l_1l_2>1)
		}{
			- \opF{\fcZ{2}\fcZ{l_3,l_4}}{2} & (l_1l_2=1 \text{\;and\ } l_3>1)
		}{
			\opF{\fcZ{3}\fcZ{l-3}}{3} - \opF{ \fcZ{2}\fcZS{1,l-3} }{2} 
			&
			(l_1l_2l_3=1 \text{\;and\ } l_4>1)
		}
	}
\end{lemma}
\envProof{
Strictly speaking,
	it is necessary to verify the fact that the polylogarithms in the lemma
	have asymptotic properties such as \refEq{2.3_Plain_EqAsymptoticForm}.
This fact is easily seen in the courses of the proofs below by virtue of \refLem[s]{2.3_LemmaA} and \ref{2.3_Lem1},
	and we do not mention it anymore.
	 
Equations \refEq{2.3_Lem2_Eq1} and \refEq{2.3_Lem2_Eq2} except the case of $l_1l_2l_3=1$ are
	proved in \cite[Lemma 2.4]{Machide12Ax1},
	and we omit their proofs.
The remaining case is derived as follows.
We see from \refEq{2.3_Lem1_Eq} that
	\envMThPd{
		\lcstC{ \lgP{1,2}{z,z^3} }
	}{
		\lcstC{ \lgP{1,2}{z,z} }
	}{
		\fcZS{1,2}
	}
By using the harmonic relation
	\envHLine
	{
		\lgP{1}{z}\lgP{1}{z}\lgP{1}{z}
	}
	{
		6 \lgP{1,1,1}{z,z^2,z^3}
		+
		3 \lgP{2,1}{z^2,z^3}
		+
		3 \lgP{1,2}{z,z^3}
		+
		\lgP{3}{z^3}
	}
	derived from the decomposition of $\Sm{m_1,m_2,m_3>0}$,
	we thus obtain
	\envOTLineCm
	{
		6 \lcstC{ \lgP{1,1,1}{z,z^2,z^3} }
	}
	{
		- 3 \fcZ{2,1} - 3\fcZS{1,2} - \fcZ{3}
	}
	{
		- 3\fcZS{1,2} - 4\fcZ{3}
	}
	where we used the simplest sum formula $\fcZ{2,1}=\fcZ{3}$. 
This proves \refEq{2.3_Lem2_Eq2} for $l_1l_2l_3=1$
	since $\fcZS{1,2}=-2\fcZ{3}$ 
	which follows from the image of the harmonic relation
	\envHLine[2.3_Lem2P_Eq1]
	{
		\lgP{1}{z}\lgP{k-1}{z}
	}
	{
		\lgP{1,k-1}{z,z^2} + \lgP{k-1,1}{z,z^2} + \lgP{k}{z^2} 
	}
	under $\lcstC[]{}$ for $k=3$.
	
Equation \refEq{2.3_Lem2_Eq3} for $l_1l_2>1$ is clear 
	because of \refEq{2.3_Plain_EqCofMPL} and \refEq{2.3_Lem1_Eq}.

We prove \refEq{2.3_Lem2_Eq3} for $l_1l_2=1$ and $l_3>1$.
By \refEq{2.1_Lem1iiP_Eq} with $z_1=z_2=z_3=z_4=z$, 
	we obtain
	\envMLinePd
	{
		\lgP{1,1}{z,z} \lgP{l_3,l_4}{z,z} 
	}
	{
		\lgP{1,1,l_3,l_4}{z,z,z^2,z^2} + \lgP{1,l_3,1,l_4}{z,z^2,z^2,z^2}
		\\
		+ 
		\lgP{1,l_3,l_4,1}{z,z^2,z^2,z^2} + \lgP{l_3,1,1,l_4}{z,z^2,z^2,z^2} + \lgP{l_3,1,l_4,1}{z,z^2,z^2,z^2}  
		\\
		+ 
		\lgP{l_3,l_4,1,1}{z,z,z^2,z^2} + \lgP{l_3+1,1,l_4}{z^2,z^2,z^2} + \lgP{l_3+1,l_4,1}{z^2,z^2,z^2} 
		\\
		+ 
		\lgP{1,l_3+1,l_4}{z,z^2,z^2} + \lgP{l_3,l_4+1,1}{z,z^2,z^2} + \lgP{1,l_3,l_4+1}{z,z^2,z^2} 
		\\
		+ \lgP{l_3,1,l_4+1}{z,z^2,z^2} + \lgP{l_3+1,l_4+1}{z^2,z^2}
	}
From  this, \refEq{2.3_Plain_EqCofMPL} and \refEq{2.3_Lem1_Eq}, it follows that
	\envMLinePd[2.3_Lem2P_Eq2]
	{
		\fcZS{1,1,l_3,l_4}
	}
	{
		-
		\bkR[b]{
			\fcZS{1,l_3,1,l_4}  + \fcZS{1,l_3,l_4,1} + \fcZS{l_3,1,1,l_4} + \fcZS{l_3,1,l_4,1} 
			\\
			+ \fcZS{l_3,l_4,1,1} + \fcZS{l_3+1,1,l_4} + \fcZS{l_3+1,l_4,1} + \fcZS{1,l_3+1,l_4} 
			\\
			+ \fcZS{l_3,l_4+1,1} + \fcZS{1,l_3,l_4+1} + \fcZS{l_3,1,l_4+1} + \fcZS{l_3+1,l_4+1}
		}
	}
On the other hand, 
	by \refEq{2.1_Lem1ii_Eq},
	we obtain
	\envMLinePd
	{
		\lgP{1,1}{\bkA{z}^2} \lgP{l_3,l_4}{\bkA{z}^2} 
	}
	{
		\lgP{1,1,l_3,l_4}{\bkA{z}^4} + \lgP{1,l_3,1,l_4}{\bkA{z}^4} + \lgP{1,l_3,l_4,1}{\bkA{z}^4} 
		\\
		+ \lgP{l_3,1,1,l_4}{\bkA{z}^4} + \lgP{l_3,1,l_4,1}{\bkA{z}^4} 
		+ \lgP{l_3,l_4,1,1}{\bkA{z}^4} + \lgP{l_3+1,1,l_4}{z^2,z^3,z^4} 
		\\
		+ \lgP{l_3+1,l_4,1}{z^2,z^3,z^4} + \lgP{1,l_3+1,l_4}{z,z^3,z^4} + \lgP{l_3,l_4+1,1}{z,z^3,z^4} 
		\\
		+ \lgP{1,l_3,l_4+1}{z,z^2,z^4} + \lgP{l_3,1,l_4+1}{z,z^2,z^4} + \lgP{l_3+1,l_4+1}{z^2,z^4}
	}
From this, \refEq{2.3_Plain_EqCofMPL}, \refEq{2.3_Lem1_Eq}, \refEq{2.3_Lem2_Eq1} and \refEq{2.3_Lem2P_Eq2},
	it also follows that
	\envHLinePd
	{
		\lcstC{ \lgP{1,1,l_3,l_4}{\bkA{z}^4} }
	}
	{
		\fcZS{1,1,l_3,l_4} - \opF{\fcZ{2}\fcZ{l_3,l_4}}{2}
	}
This proves \refEq{2.3_Lem2_Eq3} for $l_1l_2=1$ and $l_3>1$.

Noting 
	\envM{
		\fcZS{1,l-3}
	}{
		-(\fcZ{l-3,1}+\fcZ{l-2})
	}
	which is obtained by \refEq{2.3_Plain_EqCofMPL}, \refEq{2.3_Lem1_Eq} and \refEq{2.3_Lem2P_Eq1} with $k=l-2$,
	we can similarly prove \refEq{2.3_Lem2_Eq3} for $l_1l_2l_3=1$ and $l_4>1$
	by the use of \refEq{2.1_Lem1iP_Eq} with $z_1=z_2=z_3=z_4=z$
	and \refEq{2.1_Lem1i_Eq}.
We thus omit the proof.	
}

We prove \refProp{2.3_Prop1}.
\envProof[\refProp{2.3_Prop1}]{
Equations \refEq{2.3_Prop1_Eq1} and \refEq{2.3_Prop1_Eq2} are derived from 
	\refEq{2.3_Lem2_Eq1} and \refEq{2.3_Lem2_Eq2}, respectively. 
It is seen from \refEq{2.3_Lem2_Eq3} that
	\envLineThCm
	{
		\lcstC{ \gfcQPL{l}{x_1,x_2,x_3,x_4}{\bkA{z}^4} }
	}
	{
		\gfcQZVs{l}{x_1,x_2,x_3,x_4} 
		-
		\opF{\fcZ{2}}{2}\tpSm{l_3\geq2,l_4\geq1}{l_3+l_4=l-2} \fcZ{l_3,l_4} x_3^{l_3-1}x_4^{l_4-1}
		\lnAH
		+
		\bkR{ \opF{ \fcZ{3}\fcZ{l-3} }{3} -  \opF{ \fcZ{2}\fcZS{1,l-3} }{2} }x_4^{l-4}
	}
	{
		\gfcQZVs{l}{x_1,x_2,x_3,x_4} 
		-
		\opF{\fcZ{2}}{2}\tpSm{l_3,l_4\geq1}{l_3+l_4=l-2} \fcZS{l_3,l_4} x_3^{l_3-1}x_4^{l_4-1}
		+
		\opF{ \fcZ{3}\fcZ{l-3} }{3} x_4^{l-4}
	}
	which proves \refEq{2.3_Prop1_Eq3}.
}

\section{\sectThree} \label{sectThree}
We give a proof of \refThm{1_Thm1}
	by using
	\refProp[s]{2.1_Prop1}, \ref{2.2_Prop1}, \ref{2.3_Prop1} and \ref{2.3_Prop2} above,
	and 
	\refLem[s]{3_Lem1}, \ref{3_Lem2} and \ref{3_Lem3} below.
We will show the lemmas after the proof.
\envProof[\refThm{1_Thm1}]{
\refProp[s]{2.1_Prop1} and \ref{2.3_Prop1}
	yield
	\envLinePd[3_Thm1P_Eq1]
	{
		-
		\Sm{\sig\in\lgpC} \lacP{\sig}{\gfcQZVs{l}{x_1,x_2,x_3,x_4}}
		+
		\tpSm{a\geq3,b\geq1}{a+b=l} \Sm{\sig\in\lgpC} \lacPs{\sig}{ \gfcTZVs{a}{x_1,x_2,x_3} \gfcSZVs{b}{x_4} }
		\lnAH[]
		+
		\tpSm{a,b\geq2}{a+b=l} \Sm{\sig\in\lgpCb} \lacPs{\sig}{ \gfcDZVs{a}{x_1,x_2} \gfcDZVs{b}{x_3,x_4} }
		-
		\tpSm{a\geq2,b,c\geq1}{a+b+c=l} \Sm{\sig\in\lgpC} \lacPs{\sig}{ \gfcDZVs{a}{x_1,x_2} \gfcSZVs{b}{x_3} \gfcSZVs{c}{x_4} }
		\lnAH
		+
		\tpSm{a,b,c,d\geq1}{a+b+c+d=l} \gfcSZVs{a}{x_1} \gfcSZVs{b}{x_2} \gfcSZVs{c}{x_3} \gfcSZVs{d}{x_4}
		\nonumber
	}
	{
		\bkR[G]{ \tpSm{l_1,l_2,l_3,l_4\geq1}{l_1+l_2+l_3+l_4=l} x_1^{l_1-1} x_2^{l_2-1} x_3^{l_3-1} x_4^{l_4-1} } \fcZ{l}
	}
(The all extra values appearing in \refProp{2.3_Prop1} are canceled each other.)
We find from \refEq{3_Thm1P_Eq1}, \refProp[s]{2.2_Prop1} and \ref{2.3_Prop2} that
	\envLinePd
	{
		\Sm{\sig\in\lgpS} \lacP{\sig}{ \gfcQZVs{l}{x_{1234},x_{234},x_{34},x_{4}} }
		\lnAH
		-
		\Sm{\sig\in\lgpC\cup\lgpC[(34)]}
		\lacPs[g]{\sig}{
			\Sm{\rho\in\gpCyc{(234)}} \gfcQZVs{l}{x_{134},x_{\rho(2)\rho(3)\rho(4)},x_{\rho(3)\rho(4)},x_{\rho(4)}}
			\lnAHs{80}
			+
			\Sm{\rho\in\gpCyc{(24)}} \gfcQZVs{l}{x_{314},x_{14},x_{\rho(2)\rho(4)},x_{\rho(4)}}
			+
			\gfcQZVs{l}{x_{341},x_{41},x_{1},x_{2}}
		}
		\lnAH
		+
		\Sm{\sig\in\lgpC} 
		\lacPs[g]{\sig}{
			\Sm{\rho\in\gpCyc{(24)}} \gfcQZVs{l}{x_{13},x_{32},x_{\rho(2)\rho(4)},x_{\rho(4)}} 
			+
			\Sm{\rho\in\gpCyc{(34)}} \gfcQZVs{l}{x_{14},x_{24},x_{\rho(3)\rho(4)},x_{\rho(4)}}
			\lnAHs{50}
			+
			\gfcQZVs{l}{x_{13},x_{23},x_{3},x_{4}} + \gfcQZVs{l}{x_{14},x_{42},x_{2},x_{3}} + \gfcQZVs{l}{x_{41},x_{1},x_{2},x_{3}}
			\lnAHs{50}
			-
			\gfcQZVs{l}{x_{1},x_{2},x_{3},x_{4}} 
		}	
	}
	{
		\bkR[G]{ \tpSm{l_1,l_2,l_3,l_4\geq1}{l_1+l_2+l_3+l_4=l} x_1^{l_1-1} x_2^{l_2-1} x_3^{l_3-1} x_4^{l_4-1} } \fcZ{l}
	}
By \refLem[s]{3_Lem1} and \ref{3_Lem2},
	this equation holds
	if we replace $\gfcQZVs[]{l}{}$ by $\gfcQZV[]{l}{}$,
	where $\gfcQZV{l}{x_1,x_2,x_3,x_4}$ is a parameterized sum of multiple zeta values of weight $l$ which is defined by
	\envHLineDefPd[3_Thm1P_DefPSQZV]
	{
		\gfcQZV{l}{x_1,x_2,x_3,x_4}
	}
	{
		\pSmN x_1^{l_1-1}x_2^{l_2-1}x_3^{l_3-1}x_4^{l_4-1} \fcZ{l_1,l_2,l_3,l_4}
	}
We thus see from \refLem{3_Lem3} that
	\envLineCm[3_Thm1P_Eq2]
	{
		\Sm{\sig\in\lgpS} \lacP{\sig}{ \gfcQZV{l}{x_{1234},x_{234},x_{34},x_{4}} }
		\lnAH
		-
		\Sm{\sig\in\lgpC\cup\lgpC[(34)]}
		\lacPs[g]{\sig}{
			\Sm{\rho\in\gpCyc{(234)}} \gfcQZV{l}{x_{134},x_{\rho(2)\rho(3)\rho(4)},x_{\rho(3)\rho(4)},x_{\rho(4)}}
			\lnAHs{70}
			+
			\Sm{\rho\in\gpCyc{(24)}} \gfcQZV{l}{x_{314},x_{14},x_{\rho(2)\rho(4)},x_{\rho(4)}}
			+
			\gfcQZV{l}{x_{341},x_{41},x_{1},x_{2}}
		}
		\lnAH[]
		+
		\Sm{\sig\in\lgpC} 
		\lacPs[g]{\sig}{
			\Sm{\rho\in\gpCyc{(\nu(2)4)}} \Sm{\nu\in\lgpCb} \gfcQZV{l}{x_{1\nu(3)},x_{\nu(3)2},x_{\rho\nu(2)\rho(4)},x_{\rho(4)}}
			\lnAHs{40}
			+
			\Sm{\nu\in\lgpCb} \gfcQZV{l}{x_{\nu(1)3},x_{2\nu(3)},x_{\nu(3)},x_{\nu(4)}}
			+
			\gfcQZV{l}{x_{41},x_{1},x_{2},x_{3}}
			-
			\gfcQZV{l}{x_{1},x_{2},x_{3},x_{4}} 
		}	
		\nonumber
	}
	{
		\bkR[G]{ \tpSm{l_1,l_2,l_3,l_4\geq1}{l_1+l_2+l_3+l_4=l} x_1^{l_1-1} x_2^{l_2-1} x_3^{l_3-1} x_4^{l_4-1} } \fcZ{l}
	}
	which with \refEq{3_Thm1P_DefPSQZV} proves \refEq{1_Thm1_Eq1}.
}

We show the lemmas.

\begin{lemma}\label{3_Lem1}
We have
	\envHLinePd[3_Lem1_Eq]
	{
		\gfcQZVs{l}{x_1,x_2,x_3,x_4}
	}
	{
		\gfcQZV{l}{x_1,x_2,x_3,x_4} + \gfcQZVs{l}{0,x_2,x_3,x_4} 
	}
\end{lemma}
\envProof{
Equation \refEq{3_Lem1_Eq} is obvious because of \refEq{2.3_Plain_DefGfcRZV} and \refEq{3_Thm1P_DefPSQZV}.
}

\begin{lemma}\label{3_Lem2}
We have
	\envHLineCEPd
	{
		\Sm{\sig\in\lgpS} \lacP{\sig}{ \gfcQZVs{l}{0,x_{234},x_{34},x_{4}} }
	}
	{
		\Sm{\sig\in\lgpC\cup\lgpC[(34)] \atop \rho\in\gpCyc{(234)}} \lacP{\sig\rho}{ \gfcQZVs{l}{0,x_{234},x_{34},x_{4}} }
	}
	{
		\Sm{\sig\in\lgpC\cup\lgpC[(34)] \atop \rho\in\gpCyc{(24)}} \lacP{\sig}{ \gfcQZVs{l}{0,x_{14},x_{\rho(2)\rho(4)},x_{\rho(4)}} }
	}
	{
		\Sm{\sig\in\lgpC} 
		\lacPs[g]{\sig}{
			\Sm{\rho\in\gpCyc{(24)}} \gfcQZVs{l}{0,x_{32},x_{\rho(2)\rho(4)},x_{\rho(4)}} 
			\lnAHs{40}
			+
			\Sm{\rho\in\gpCyc{(34)}} \gfcQZVs{l}{0,x_{24},x_{\rho(3)\rho(4)},x_{\rho(4)}}
		}
	}
	{
		\Sm{\sig\in\lgpC\cup\lgpC[(34)]} \lacP{\sig}{ \gfcQZVs{l}{0,x_{41},x_{1},x_{2}} }
	}
	{
		\Sm{\sig\in\lgpC} \lacPs[B]{\sig}{ \gfcQZVs{l}{0,x_{23},x_{3},x_{4}} + \gfcQZVs{l}{0,x_{42},x_{2},x_{3}} }
	}
	{
		\Sm{\sig\in\lgpC} \lacP{\sig}{ \gfcQZVs{l}{0,x_{1},x_{2},x_{3}} }
	}
	{
		\Sm{\sig\in\lgpC} \lacP{\sig}{ \gfcQZVs{l}{0,x_{2},x_{3},x_{4}} }
	}
\end{lemma}
\envProof{
We use \refEq{2.1_Lem2_Eq1} without notice.
We obtain the first equation by
	\envMThPd
	{
		(\lgpC\cup\lgpC[(34)]) \cdot \gpCyc{(234)}
	}
	{
		\lgpC\cup\lgpC[(234)]\cup\lgpC[(243)]\cup\lgpC[(34)]\cup\lgpC[(24)]\cup\lgpC[(23)]
	}
	{
		\lgpS
	}
Since $\lgpC=\lgpC[(13)(24)]$ and $\lgpC[(34)]=\lgpC[(123)]$,
	we have
	\envHLineF
	{
		\Sm{\sig\in\lgpC \atop \rho\in\gpCyc{(24)}} \lacP{\sig}{ \gfcQZVs{l}{0,x_{14},x_{\rho(2)\rho(4)},x_{\rho(4)}} }
	}
	{
		\Sm{\sig\in\lgpC[(13)(24)] \atop \rho\in\gpCyc{(24)}} \lacP{\sig}{ \gfcQZVs{l}{0,x_{14},x_{\rho(2)\rho(4)},x_{\rho(4)}} }
	}
	{
		\Sm{\sig\in\lgpC \atop \rho\in(13)(24)\cdot\gpCyc{(24)}} \lacP{\sig}{ \gfcQZVs{l}{0,x_{32},x_{\rho(2)\rho(4)},x_{\rho(4)}} }
	}
	{
		\Sm{\sig\in\lgpC \atop \rho\in\gpCyc{(24)}} \lacP{\sig}{ \gfcQZVs{l}{0,x_{32},x_{\rho(2)\rho(4)},x_{\rho(4)}} }
	}
	and
	\envHLineFCm
	{
		\Sm{\sig\in\lgpC[(34)] \atop \rho\in\gpCyc{(24)}} \lacP{\sig}{ \gfcQZVs{l}{0,x_{14},x_{\rho(2)\rho(4)},x_{\rho(4)}} }
	}
	{
		\Sm{\sig\in\lgpC[(123)] \atop \rho\in\gpCyc{(24)}} \lacP{\sig}{ \gfcQZVs{l}{0,x_{14},x_{\rho(2)\rho(4)},x_{\rho(4)}} }
	}
	{
		\Sm{\sig\in\lgpC \atop \rho\in(123)\cdot\gpCyc{(24)}} \lacP{\sig}{ \gfcQZVs{l}{0,x_{24},x_{\rho(2)\rho(4)},x_{\rho(4)}} }
	}
	{
		\Sm{\sig\in\lgpC \atop \rho\in\gpCyc{(34)}}  \lacP{\sig}{ \gfcQZVs{l}{0,x_{24},x_{\rho(3)\rho(4)},x_{\rho(4)}} }
	}
	which gives the second one.
Similarly the third one follows.
The fourth one is derived from $\lgpC=\lgpC[(1234)]$.
}
\begin{lemma}\label{3_Lem3}
We have
	\envMLine[3_Lem3_Eq1]
	{
		\Sm{\sig\in\lgpC} 
		\lacPs[g]{\sig}{
			\Sm{\rho\in\gpCyc{(24)}} \gfcQZV{l}{x_{13},x_{32},x_{\rho(2)\rho(4)},x_{\rho(4)}} 
			+
			\Sm{\rho\in\gpCyc{(34)}} \gfcQZV{l}{x_{14},x_{24},x_{\rho(3)\rho(4)},x_{\rho(4)}}
		}
		\\
	}
	{
		\Sm{\sig\in\lgpC} \lacPs[g]{\sig}{ \Sm{\rho\in\gpCyc{(\nu(2)4)}} \Sm{\nu\in\lgpCb} \gfcQZV{l}{x_{1\nu(3)},x_{\nu(3)2},x_{\rho\nu(2)\rho(4)},x_{\rho(4)}} }
	}
	and
	\envMLinePd[3_Lem3_Eq2]
	{
		\Sm{\sig\in\lgpC} \lacPs{\sig}{ \gfcQZV{l}{x_{13},x_{23},x_{3},x_{4}} + \gfcQZV{l}{x_{14},x_{42},x_{2},x_{3}} }
		\\
	}
	{
		\Sm{\sig\in\lgpC} \lacPs[g]{\sig}{ \Sm{\nu\in\lgpCb} \gfcQZV{l}{x_{\nu(1)3},x_{2\nu(3)},x_{\nu(3)},x_{\nu(4)}} }
	}
\end{lemma}
\envProof{
Since $\lgpCb=\bkB{\gpu, (1234)}$,
	we have
	\envLineCm
	{
		\Sm{\rho\in\gpCyc{(24)}} \gfcQZV{l}{x_{13},x_{32},x_{\rho(2)\rho(4)},x_{\rho(4)}} 
		+
		\Sm{\rho\in\gpCyc{(34)}} \gfcQZV{l}{x_{14},x_{24},x_{\rho(3)\rho(4)},x_{\rho(4)}}
	}
	{
		\Sm{\rho\in\gpCyc{(\nu(2)4)}} \Sm{\nu\in\lgpCb} \gfcQZV{l}{x_{1\nu(3)},x_{\nu(3)2},x_{\rho\nu(2)\rho(4)},x_{\rho(4)}} 
	}
	which gives \refEq{3_Lem3_Eq1}.
We also have
	\envHLinePd
	{
		\gfcQZV{l}{x_{13},x_{23},x_{3},x_{4}}
		+
		\gfcQZV{l}{x_{32},x_{24},x_{4},x_{1}}
	}
	{
		\Sm{\nu\in\lgpCb} \gfcQZV{l}{x_{\nu(1)3},x_{2\nu(3)},x_{\nu(3)},x_{\nu(4)}}
	}
This proves \refEq{3_Lem3_Eq2}
	because of
	\envHLine
	{
		\Sm{\sig\in\lgpC} \lacP{\sig}{ \gfcQZV{l}{x_{32},x_{24},x_{4},x_{1}} }
	}
	{
		\Sm{\sig\in\lgpC} \lacP{\sig}{ \gfcQZV{l}{x_{14},x_{42},x_{2},x_{3}} }
	}
	which follows from $\lgpC=\lgpC[(13)(24)]$.
}

\section{\sectFour} \label{sectFour}
In this final section, 
	we derive \refThm{1_Thm2} from \refThm{1_Thm1}.
Before proving \refThm{1_Thm2},
	we prepare some equations by substituting $0$ or $1$ for each parameter $x_j$ in \refEq{1_Thm1_Eq1}.
For the substitutions,
	the mathematical software ``Maxima'' is used implicitly,
	and \refEq{3_Thm1P_Eq2} instead of \refEq{1_Thm1_Eq1} is referred to 
	since both are equivalent and \refEq{3_Thm1P_Eq2} is convenient to calculate.
\begin{lemma}\label{4_Lem1}
Let $l$ be a positive integer with $l\geq5$.
We have
	\envHLineCENmePd
	{\label{4_Lem1_Eq1}
		\envSLine{
			2\gfcQZV{l}{2,2,2,1} + \gfcQZV{l}{2,2,1,1} + \gfcQZV{l}{2,1,1,1} 
			\\
			- 
			2\gfcQZV{l}{1,2,2,1} - \gfcQZV{l}{1,2,1,1} - 3\gfcQZV{l}{1,1,1,1}
		}
	}
	{
		(l-3)\fcZ{l}
	}
	{\label{4_Lem1_Eq2}
		\envSLine{
			4\gfcQZV{l}{2,2,2,1} + 2\gfcQZV{l}{2,2,1,1} - 4\gfcQZV{l}{1,2,2,1} 
			\\
			- 2\gfcQZV{l}{1,2,1,1} - 4\gfcQZV{l}{1,1,2,1}
		}
	}
	{
		(l-3)\fcZ{l}
	}
	{\label{4_Lem1_Eq3}
		\envSLine{
			6\gfcQZV{l}{3,3,2,1} + 4\gfcQZV{l}{3,2,2,1} + 2\gfcQZV{l}{3,2,1,1} 
			\\
			-
			6\gfcQZV{l}{2,3,2,1} -8\gfcQZV{l}{2,2,2,1} - 4\gfcQZV{l}{2,2,1,1} 
			\\
			- 
			2\gfcQZV{l}{2,1,2,1} - 2\gfcQZV{l}{2,1,1,1} + 4\gfcQZV{l}{1,2,2,1} 
			\\
			+ 
			2\gfcQZV{l}{1,2,1,1} + 2\gfcQZV{l}{1,1,2,1} + 3\gfcQZV{l}{1,1,1,1}
		}
	}
	{
		\binom{l-2}{2} \fcZ{l}
	}
	{\label{4_Lem1_Eq4}
		\envSLine{
			24\gfcQZV{l}{4,3,2,1} - 24\gfcQZV{l}{3,3,2,1} - 16\gfcQZV{l}{3,2,2,1} 
			\\
			- 
			8\gfcQZV{l}{3,2,1,1} + 16\gfcQZV{l}{2,2,2,1} + 8\gfcQZV{l}{2,2,1,1} 
			\\
			+ 
			4\gfcQZV{l}{2,1,1,1} - 4\gfcQZV{l}{1,1,1,1}
		}
	}
	{
		\binom{l-1}{3}\fcZ{l}
	}
\end{lemma}
\envProof{
Equations \refEq{4_Lem1_Eq1}, \refEq{4_Lem1_Eq2}, \refEq{4_Lem1_Eq3} and \refEq{4_Lem1_Eq4}
	are obtained by substituting 
	$(1,1,0,0)$, $(1,0,1,0)$, $(1,1,1,0)$ and $(1,1,1,1)$ for $(x_1,x_2,x_3,x_4)$ in \refEq{3_Thm1P_Eq2}, respectively.
}

We prove \refThm{1_Thm2}.
Note that
	\envM{
		\gfcQZV{l}{1,1,1,1}
	}{
		\fcZ{l}
	}
	by \refEq{1_Pl_EqSF} with $n=4$.
\envProof[\refThm{1_Thm2}]{
Formulas \refEq{1_Thm2i_Eq1} and \refEq{1_Thm2i_Eq2} respectively follow from \refEq{4_Lem1_Eq1} and \refEq{4_Lem1_Eq2}
	because of \refEq{3_Thm1P_DefPSQZV}.
Formula \refEq{1_Thm2i_Eq3}
	is derived from subtracting 
	\refEq{1_Thm2i_Eq2} from twice \refEq{1_Thm2i_Eq1}.

We prove \refEq{1_Thm2ii_Eq1} next.
Since
	\envMCm{
		\opF[s]{3(l-3)}{2} + \binom{l-2}{2}
	}{
		\opF[s]{(l+1)(l-3)}{2}
	}
	adding up \refEq{4_Lem1_Eq1}, the half of \refEq{4_Lem1_Eq2}, and \refEq{4_Lem1_Eq3} yields
	\envMLinePd[4_Thm2P_Eq]
	{
		6\gfcQZV{l}{3,3,2,1} + 4\gfcQZV{l}{3,2,2,1} + 2\gfcQZV{l}{3,2,1,1} - 6\gfcQZV{l}{2,3,2,1} - 4\gfcQZV{l}{2,2,2,1} 
		\\
		- 
		2\gfcQZV{l}{2,2,1,1} - 2\gfcQZV{l}{2,1,2,1} - \gfcQZV{l}{2,1,1,1}
	}
	{
		\opF{(l+1)(l-3)}{2} \fcZ{l}
	}
Since 
	\envMCm{
		\opF[s]{(l+1)(l-3)}{2}  + \opF[s]{ \binom{l-1}{3} }{4} + 1
	}{
		\opF[s]{ (l+1)(l^2+5l-18) }{24}
	}
	adding up \refEq{4_Thm2P_Eq} and the quarter of \refEq{4_Lem1_Eq4} also yields
	\envHLineCm
	{
		6\gfcQZV{l}{4,3,2,1} - 6\gfcQZV{l}{2,3,2,1} - 2\gfcQZV{l}{2,1,2,1}
	}
	{
		\opF{ (l+1)(l^2+5l-18) }{24} \fcZ{l}
	}
	which proves \refEq{1_Thm2ii_Eq1}.
}
\begin{remark}\label{4_Rem1}
We rewrite \refEq{4_Lem1_Eq3}, \refEq{4_Lem1_Eq4} and \refEq{4_Thm2P_Eq},
	which are necessary to prove \refEq{1_Thm2ii_Eq1},
	in terms of quadruple zeta values $\fcZ{\lbfL}=\fcZ{l_1,l_2,l_3,l_4}$
	to see the explicit relations among these values.
	\envHLineCSCmNme
	{\label{4_Rem1_Eq1}
		\envSLine{
			\pSmN
			\bkR[b]{
				3^{l_{12}-1}2^{l_3} +  3^{l_1-1}2^{l_{23}} + 3^{l_1-1}2^{l_2} 
				\slnAHs{25}
				- 
				3^{l_2} 2^{l_{13}-1} - 2^{l_{123}} - 2^{l_{12}} - 2^{l_{13}-1} - 2^{l_1} 
				\slnAHs{25}
				+
				2^{l_{23}} + 2^{l_2} + 2^{l_3}
			}
			\fcZ{\lbfL}
		}
	}
	{
		\opF{l(l-5)}{2}\fcZ{l}
	}
	{\label{4_Rem1_Eq2}
		\envSLine{
			\pSmN
			\bkR[b]{
				3^{l_2}2^{2l_1+l_3-1} - 3^{l_{12}-1}2^{l_3+1} - 3^{l_1-1}2^{l_{23}+1}
				\slnAHs{25}
				- 
				3^{l_1-1}2^{l_2+1} + 2^{l_{123}} + 2^{l_{12}} + 2^{l_1}
			}
			\fcZ{\lbfL}
		}
	}
	{
		\opF{(l+1)(l^2-7l+18)}{12}\fcZ{l}
	}
	{\label{4_Rem1_Eq3}
		\envSLine{
			\pSmN
			\bkR[b]{
				3^{l_{12}-1}2^{l_3+1} + 3^{l_1-1}2^{l_{23}+1} + 3^{l_1-1}2^{l_2+1} 
				\slnAHs{25}
				- 
				3^{l_2}2^{l_{13}} - 2^{l_{123}} - 2^{l_{12}} - 2^{l_{13}} - 2^{l_1} 
			}
			\fcZ{\lbfL}
		}
	}
	{
		(l+1)(l-3) \fcZ{l}
	}
	where 
	\refEq{4_Rem1_Eq1}, \refEq{4_Rem1_Eq2} and \refEq{4_Rem1_Eq3}  
	correspond to \refEq{4_Lem1_Eq3}, \refEq{4_Lem1_Eq4} and \refEq{4_Thm2P_Eq}, respectively.
\end{remark}

\section*{Acknowledgements}
The author would like to thank the National Center for Theoretical Sciences (Taiwan) 
	for the hospitality and support
	in his visit for ``2012 NCTS Japan-Taiwan Joint Conference on Number Theory''.
The author is also thankful
	to Professor Wen-Chin Liaw
	for introducing him the preprint \cite{OEL12Ax}.


\begin{flushleft}
\mbox{}\\ \qquad \mbox{}\\ \qquad
\majorName		\mbox{}\\ \qquad
\departmentName	\mbox{}\\ \qquad
\organizationName	\mbox{}\\ \qquad
\placeAddress		\mbox{}\\ \qquad
\emailAddress
\end{flushleft}

\end{document}